\documentclass[12pt,reqno]{amsart}
\usepackage{a4wide}
\usepackage{amsfonts,amsmath,amssymb,amscd}
\usepackage{amsthm}
\usepackage{amsmath,todonotes}
\usepackage{amscd}
\usepackage{cancel}
\usepackage{verbatim}
\usepackage{color}
\usepackage[all]{xy}
\usepackage[mathscr]{eucal}
\usepackage{mathrsfs}
\usepackage{fullpage}
\usepackage{enumitem}
\usepackage{hyperref}
\usepackage{bbm}
\usepackage{qtree}
\usepackage{mathrsfs}
\usepackage{bm}
\usepackage{bigints}
\usepackage{url}
\allowdisplaybreaks

\newcommand{\pFq}[5]{\ensuremath{{}_{#1}F_{#2} \left( \genfrac{}{}{0pt}{}{#3}
		{#4} \bigg| {#5} \right)}}
\renewcommand{\a}{\alpha}
\renewcommand{\b}{\beta}

\renewcommand{\l}{\lambda}

\usepackage[headheight=110pt,top=0.6in, bottom=0.7in, left=0.6in, right=0.6in]{geometry}
\numberwithin{equation}{section}

\newtheorem{theorem}{Theorem}

\theoremstyle{definition}

\numberwithin{theorem}{section} 
\numberwithin{equation}{section}
\numberwithin{table}{section}

\newcommand{\g}{\gamma}

\allowdisplaybreaks

\makeatletter
\def\proof{\@ifnextchar[{\@oproof}{\@nproof}}
\def\@oproof[#1][#2]{\trivlist\item[\hskip\labelsep\textit{#2 Proof of\
		#1.}~]\ignorespaces}
\def\@nproof{\trivlist\item[\hskip\labelsep\textit{Proof.}~]\ignorespaces}
\begin{document}
	\title[]{Recent developments pertaining to Ramanujan's formula for odd zeta values}
		
	\author{Atul Dixit}
	\address{Department of Mathematics, Indian Institute of Technology Gandhinagar, Palaj, Gandhinagar 382355, Gujarat, India} 
	\email{adixit@iitgn.ac.in}
	
	\thanks{2020 \textit{Mathematics Subject Classification.} Primary 11M06, 44A20; Secondary 11J91\\
		\textit{Keywords and phrases.} Lambert series, odd zeta values, Eisenstein series, Ramanujan polynomials}

	\medskip
	\begin{abstract}
	In this expository article, we discuss the contributions made by several mathematicians with regard to a famous formula of Ramanujan for odd zeta values. The goal is to complement the excellent survey by Berndt and Straub \cite{berndtstraubzeta} with some of the recent developments that have taken place in the area in the last decade or so.
	\end{abstract}
	
	\dedicatory{Dedicated to Srinivasa Ramanujan on the occasion of his $136^{\textup{th}}$ birthday}
	\maketitle	
	
\tableofcontents
\section{Introduction}
At the Ramanujan centenary conference at the University of Illinois at Urbana-Champaign in 1987, the renowned physicist and mathematician Freeman Dyson remarked, \emph{``That was the wonderful thing about Ramanujan. He discovered so much, and yet he left so much more in his garden for other people to discover''.} Dyson's remarks hold true to this day as can be seen from the stunning developments currently happening in the theory of mock modular forms (having its origin in Ramanujan's mock theta functions), the theory of partitions, Rogers-Ramanujan identities (hard hexagon model, Kanade-Russell conjectures), to name a few. Another result of Ramanujan which is in the same league as above is his famous formula for the values of the Riemann zeta function at odd integers (other than $1$), which we will shortly describe. 

The Riemann zeta function $\zeta(s)$ is defined for Re$(s)>1$ by the absolutely convergent series
\begin{equation*}
\zeta(s):=\sum_{n=1}^{\infty}\frac{1}{n^s}.
\end{equation*}
It converges uniformly in Re$(s)\geq1+\epsilon$ for any $\epsilon>0$ and hence represents an analytic function in Re$(s)>1$	. Further, it is well-known that $\zeta(s)$ can be analytically continued in the entire complex plane except for a simple pole at $s=1$. The Riemann zeta function is just a simple case of what is known as an $L$-function, which, in the half-plane Re$(s)>\sigma_a$, is defined by 
\begin{equation*}
	\sum_{n=1}^{\infty}\frac{a(n)}{n^s},
\end{equation*}
where $a(n)$ is an arithmetic function and $\sigma_a$ is the abscissa of absolute convergence of the series.

For $m\geq1$, it is well-known that $\zeta(-2m)=0$. Euler gave an explicit representation for $\zeta(2m)$, $m\geq1$, given by
\begin{equation}\label{zetaevenint}
	\zeta(2 m ) = (-1)^{m +1} \frac{(2\pi)^{2 m}B_{2 m }}{2 (2 m)!},
\end{equation}
where $B_{m}$ is the $m^{\textup{th}}$ Bernoulli number defined by
\begin{equation*}
	\frac{x}{e^{x}-1}=:\sum_{m=1}^{\infty}\frac{B_m}{m!}x^m,\hspace{5mm}(|x|<2\pi),
\end{equation*} 
This implies, in particular, $\zeta(2)=\pi^2/6, \zeta(4)=\pi^4/90, \zeta(6)=\pi^6/945$ etc.
One of the reasons why Euler's formula is extremely important is that it at once implies that $\zeta(2m)$ for any $m\geq 1$ is transcendental. This is owing to the fact that $\pi$ is transcendental and $B_m, m\geq1,$ is a rational number. In particular, this implies  that  for any $m\geq 1$, $\zeta(2m)$ is irrational.

Now a natural question arises here - \emph{what can we say about $\zeta(2m+1)$ for any $m\geq 1$?}\footnote{It is well-known for $m\geq0$ that $\zeta(-2m-1)=-\frac{B_{2m+2}}{2m+2}$, a rational number.} (We shall henceforth call $\zeta(2m+1)$ as the odd zeta values.) Except in the case of $\zeta(3)$, the answer is far from being known! Even in the case of $\zeta(3)$, we only know, thanks to Ap\'{e}ry \cite{apery1}, \cite{apery2}, that $\zeta(3)$ is irrational. Nothing is known currently on its possibility of being transcendental. The reader is encouraged to see the paper of Rajkumar \cite{rajkumar} and that of Zudilin \cite{zudilin1} for interesting accounts on how a result of Ramanujan inspired Ap\'{e}ry in his proof of irrationality of $\zeta(3)$. As far as the other odd zeta values $\zeta(5), \zeta (7), \zeta(9),\cdots$ are concerned, we do not even know if they are irrational. However, it is known due to Zudilin \cite{zudilin} that at least one of  $\zeta(5), \zeta (7), \zeta(9)$ or $\zeta(11)$ must be irrational. Moreover, Rivoal \cite{rivoal} has shown that infinitely many numbers in the set $\{\zeta(2m+1)\}_{m=2}^{\infty}$ are irrational, but we do not know exactly which ones. 

The topic of evaluating the Riemann zeta function, and more generally $L$-functions, at special values of their arguments remains to this day an evergreen topic. While there may not always be exact formulas that are available, there are transformations in the literature involving such functions which are fundamentally important.

One such is the formula of Ramanujan for $\zeta(2m+1)$ alluded to in the second paragraph of the introduction. Let $\a, \b>0$ with $\a\b=\pi^2$ and $m\in\mathbb{Z}\backslash\{0\}$. Then Ramanujan's formula is given by\footnote{Ramanujan's formula actually holds for any complex $\a, \b$ such that $\textup{Re}(\a)>0, \textup{Re}(\b)>0$ and $\a\b=\pi^2$.} \cite[p.~173, Ch. 14, Entry 21(i)]{ramnote}, \cite[p.~319-320, formula (28)]{lnb}, \cite[p.~275-276]{bcbramsecnote} 
\begin{align}\label{rameqn}
	\alpha^{-m}\left\{\frac{1}{2}\zeta(2m+1)+\sum_{n=1}^\infty \frac{n^{-2m-1}}{e^{2n\alpha}-1}\right\}&=(-\beta)^{-m}\left\{\frac{1}{2}\zeta(2m+1)+\sum_{n=1}^\infty\frac{n^{-2m-1}}{e^{2n\beta}-1}\right\}\nonumber\\
	&\qquad-2^{2m}\sum_{k=0}^{m+1}\frac{(-1)^{k}B_{2k}B_{2m+2-2k}}{(2k)!(2m+2-2k)!}\alpha^{m+1-k}\beta^k.
\end{align}
The series occurring in \eqref{rameqn} are examples of Lambert series. A \emph{Lambert series} associated with the arithmetic function $a(n)$ is defined by 	$\displaystyle\sum_{n=1}^{\infty}a(n)\frac{q^{n}}{1-q^{n}}$, where $|q|<1$. If we let  $q=e^{-y}$, where Re$(y)>0$, then the above Lambert series can be written in the equivalent form
$\displaystyle\sum_{n=1}^{\infty}\frac{a(n)}{e^{ny}-1}$, which is what appears in \eqref{rameqn} with $a(n)=n^{-2m-1}$. 

There are several applications of Ramanujan's formula. We begin with its special case - a formula due to Lerch \cite{lerch}.  For  odd $m\in\mathbb{N}$, it is given by
\begin{align}\label{lerch}
	\zeta(2m+1)+2\sum_{n=1}^{\infty}\frac{n^{-2m-1}}{e^{2\pi n}-1}=\pi^{2m+1}2^{2m}\sum_{j=0}^{m+1}\frac{(-1)^{j+1}B_{2j}B_{2m+2-2j}}{(2j)!(2m+2-2j)!}.
\end{align} 
It can be thought of as a formula for evaluating special values of the Dirichlet series associated with $\textup{coth}(\pi n)$; see, for example, the paper of Straub \cite{straub}, where such evaluations are taken up further.

While the right-hand side of \eqref{lerch} is again a transcendental number, the left-hand side, unlike \eqref{zetaevenint}, is not just $\zeta(2m+1)$ but instead 
the sum of $\zeta(2m+1)$ and a rapidly convergent series. Hence, at best, we can only say \cite{gmr} that at least one of $\zeta(2m+1)$ or $\sum_{n=1}^{\infty}\frac{n^{-2m-1}}{e^{2\pi n}-1}$ is transcendental.

Ramanujan's formula encodes the fundamental transformation properties of the Eisenstein series on the full modular group SL$_2(\mathbb{Z})$. The Eisenstein series of even integral weight $k\geq2$ over SL$_2(\mathbb{Z})$ have the following Fourier series expansion:
\begin{equation}
	E_k(z)= 1- \frac{2k}{B_k}\sum_{n=1}^\infty \sigma_{k-1}(n) e^{2\pi i nz},\quad z\in\mathbb{H}\hspace{1mm}\text{(the upper half-plane)}, \nonumber
\end{equation}
where $\displaystyle\sigma_s(n)=\sum_{d|n}d^s$ is the generalized divisor function and $B_k$ are the Bernoulli numbers. Now letting $m=-\ell$ in \eqref{rameqn} gives, for $\ell>1$, the transformation formula satisfied by the Eisenstein series 
\begin{align}\label{es}
	\alpha^\ell\sum_{n=1}^\infty\frac{n^{2\ell-1}}{e^{2n\alpha}-1}-(-\beta)^\ell\sum_{n=1}^\infty\frac{n^{2\ell-1}}{e^{2n\beta}-1}=\left(\alpha^\ell-(-\beta)^\ell\right)\frac{B_{2\ell}}{4\ell}
\end{align}
Indeed, this is a reformulation of the well-known modular relation valid for $\ell>1$:
\begin{align*}
	E_{2\ell}(-1/z)=z^{2\ell}E_{2\ell}(z),\quad z\in\mathbb{H}.
\end{align*}
Moreover, when $m=-1$ in \eqref{rameqn}, we get an equivalent version of the modular transformation $E_2(-1/z)=z^{2}E_2(z)+\frac{6z}{\pi i}$ satisfied by the quasi-modular form $E_2(z)$, namely, \begin{align*}
	\a\sum_{n=1}^{\infty}\frac{n}{e^{2n\a}-1}+\b \sum_{n=1}^{\infty}\frac{n}{e^{2n\b}-1}=\frac{\a+\b}{24}-\frac{1}{4}.
\end{align*}
There are further important corollaries resulting from Ramanujan's formula. One such concerns the Eichler integral associated with the Eisenstein series $E_k(z)$. An Eichler integral corresponding to $E_k(z)$ is defined as the $(k-1)^{\textup{st}}$ primitive of $E_k(z)$. Then Ramanujan's formula for positive integers $m$ gives the transformation of Eichler integrals corresponding to $E_{2m+2}(z)$. Moreover, even though the pole of $\zeta(2m+1)$ at $m=0$ does not permit us letting  $m=0$ in \eqref{rameqn}, the Lambert series occurring in \eqref{rameqn} make perfect sense. Thus, starting with one of these Lambert series with $m=0$ and converting it into a line integral using the Perron inversion formula followed by appropriately shifting the line of integration, using the Cauchy residue theorem and taking into account the contributions of the residues at the poles, leads to
\begin{align}\label{etatransequiv}
	\sum_{n=1}^{\infty}\frac{1}{n(e^{2n\a}-1)}-\sum_{n=1}^{\infty}\frac{1}{n(e^{2n\b}-1)}&=\frac{\b-\a}{12}+\frac{1}{4}\log\left(\frac{\a}{\b}\right).
\end{align}
This transformation is nothing but an equivalent form \footnote{To see this, take logarithm on both sides of \eqref{etatrans}, use the Taylor expansion of logarithm, then substitute $\alpha=-\pi iz$ and $\beta=\pi i/z$ (so that Re$(\alpha)>0$ and Re$(\beta)>0$, and simplify.} of
\begin{align}\label{etatrans}
\eta(-1/z)=\sqrt{-iz}\eta(z),
\end{align}
 where $\eta(z)$ is the Dedekind eta-function defined for $z\in\mathbb{H}$ by $\eta(z)=e^{\frac{\pi iz}{12}}\prod_{n=1}^{\infty}(1-e^{2\pi inz})$. It is well-known that $\eta(z)$ is a half-integral weight modular form on SL$_2(\mathbb{Z})$ with its transformations twisted by roots of unity.

Thus, Ramanujan's formula \eqref{rameqn} encapsulates the  transformations of Eisenstein series on SL$_2(\mathbb{Z})$ as well as those of their Eichler integrals, and, in addition, also gives the Dedekind eta-function transformation. In fact, Ramanujan's formula is equivalent to a certain representation for period polynomials of Eisenstein series of even weight on $\textup{SL}_{2}\left(\mathbb{Z}\right)$; see \cite[Proposition 5.2, Equation (35)]{berndtstraubmathz}.

Apart from the applications in modular forms, Kirschenhofer and Prodinger \cite{kirprod} have given an application of Ramanujan's formula in theoretical computer science, in particular, in the analysis of special data structures and algorithms. To be more specific,  \eqref{rameqn} and its aforementioned corollaries are used to achieve certain distribution results on random variables related to dynamic data structures called `tries'.

The literature on Ramanujan's formula and its very many generalizations and analogues is vast. There have also been many surveys and expository articles written on the subject, for example, those by Berndt \cite{berndt-jupiter}, Berndt and Straub \cite{berndtstraubzeta} and by Zudilin \cite{zudilin1}. In this survey, we will be concentrating on the developments that have taken place in the area in the last decade or so, in particular, apr\`{e}s the excellent survey \cite{berndtstraubzeta}. To keep the survey short, we do not duplicate the material in \cite{berndtstraubzeta} concerning how Ramanujan may have proved \eqref{rameqn}. The contributions by various mathematicians  concerning \eqref{rameqn} before the last decade have been well-documented in \cite[p.~276]{bcbramsecnote}, \cite{berndtstraubzeta} and hence, barring a few exceptions, do not form a topic of discussion in this paper. 

\section{Generalizations of Ramanujan's formula \eqref{rameqn}}

The generalizations of Ramanujan's formula given here are not in chronological order.\\

\subsection{A generalization associated with the Lambert series of $n^{s}$ and its applications}
Recently, Kesarwani, Kumar and the current author \cite{dkk} considered the generalized Lambert series
	\begin{equation}\label{rames}
	\sum_{n=1}^{\infty}\frac{n^s}{e^{ny}-1}=\sum_{n=1}^{\infty}\sigma_s(n)e^{-ny},
\end{equation}
where $s\in\mathbb{C}$, Re$(y)>0$ and, for Re$(s)>-1$, obtained the following generalization \cite[Theorem 2.4]{dkk} of Ramanujan's \eqref{rameqn}.
\begin{theorem}\label{befac}
\begin{align}\label{maineqn}
	&\sum_{n=1}^\infty  \sigma_s(n)e^{-ny}+\frac{1}{2}\left(\left(\frac{2\pi}{y}\right)^{1+s}\mathrm{cosec}\left(\frac{\pi s}{2}\right)+1\right)\zeta(-s)-\frac{1}{y}\zeta(1-s)\nonumber\\
	&=\frac{2\pi}{y\sin\left(\frac{\pi s}{2}\right)}\sum_{n=1}^\infty \sigma_{s}(n)\Bigg(\frac{(2\pi n)^{-s}}{\Gamma(1-s)} {}_1F_2\left(1;\frac{1-s}{2},1-\frac{s}{2};\frac{4\pi^4n^2}{y^2} \right) -\left(\frac{2\pi}{y}\right)^{s}\cosh\left(\frac{4\pi^2n}{y}\right)\Bigg),
\end{align}
where ${}_1F_2(a;b,c;z):=\sum_{n=0}^\infty\frac{(a)_nz^n}{(b)_n(c)_nn!}$ with $z\in\mathbb{C}$ is the generalized hypergeometric function with $\ (a)_n=\frac{\Gamma(a+n)}{\Gamma(a)}$, and $\Gamma(z)$ is the Euler Gamma function. 
\end{theorem}
This generalization can be considered as a transformation for a ``complex'' analogue of Eisenstein series. They used analytic continuation to extended the validity of the above formula in a larger region and obtained the following result \cite[Theorem 2.5]{dkk}.
\begin{theorem}\label{extendedid}
	Let $m\in\mathbb{N}\cup\{0\}$ and $\textup{Re}(y)>0$. Then for $\mathrm{Re}(s)>-2m-3$, we have
	{\allowdisplaybreaks\begin{align}\label{extendedideqn}
			&\sum_{n=1}^\infty  \sigma_s(n)e^{-ny}+\frac{1}{2}\left(\left(\frac{2\pi}{y}\right)^{1+s}\mathrm{cosec}\left(\frac{\pi s}{2}\right)+1\right)\zeta(-s)-\frac{\zeta(1-s)}{y}\nonumber\\
			&=\frac{2\sqrt{2\pi}}{y^{1+\frac{s}{2}}}\sum_{n=1}^\infty\sigma_s(n)n^{-\frac{s}{2}}\left\{{}_{\frac{1}{2}}{K}_{\frac{s}{2}}\left(\frac{4\pi^2n}{y},0\right)-\frac{\pi2^{\frac{3}{2}+s}}{\sin\left(\frac{\pi s}{2}\right)}\left(\frac{4\pi^2n}{y}\right)^{-\frac{s}{2}-2}A_m\left(\frac{1}{2},\frac{s}{2},0;\frac{4\pi^2n}{y}\right)\right\}\nonumber\\
			&\qquad-\frac{y(2\pi)^{-s-3}}{\sin\left(\frac{\pi s}{2}\right)}\sum_{k=0}^m\frac{\zeta(s+2k+2)\zeta(2k+2)}{\Gamma(-s-1-2k)}\left(\frac{4\pi^2}{y}\right)^{-2k},
	\end{align}}
	where ${}_{\mu}K_{\nu}(z, w)$ is a generalized modified Bessel function defined for $\nu\in\mathbb{C}\backslash\left(\mathbb{Z}\backslash\{0\}\right)$, and $z, \mu, w\in\mathbb{C}$ such that $\mu+w\neq-\frac{1}{2}, -\frac{3}{2}, -\frac{5}{2},\cdots$, by
	%\backslash\{\cdots, -3, -2, -1, 1, 2, 3, \cdots\}
	\begin{align}\label{def2varbessel}
		{}_{\mu}K_{\nu}(z, w)&:=\frac{\pi z^w 2^{\mu+\nu-1}}{\sin(\nu\pi)}\bigg\{\left(\frac{z}{2}\right)^{-\nu}\frac{\Gamma(\mu+w+\tfrac{1}{2})}{\Gamma(1-\nu)\Gamma(w+\tfrac{1}{2}-\nu)}\pFq12{\mu+w+\tfrac{1}{2}}{w+\tfrac{1}{2}-\nu,1-\nu}{\frac{z^2}{4}}\nonumber\\
		&\quad\quad\quad\quad\quad\quad-\left(\frac{z}{2}\right)^{\nu}\frac{\Gamma(\mu+\nu+w+\tfrac{1}{2})}{\Gamma(1+\nu)\Gamma(w+\tfrac{1}{2})}\pFq12{\mu+\nu+w+\tfrac{1}{2}}{w+\tfrac{1}{2},1+\nu}{\frac{z^2}{4}}\bigg\},
	\end{align}
	with ${}_{\mu}K_{0}(z, w)=\lim_{\nu\to0}{}_{\mu}K_{\nu}(z, w)$, and
	\begin{align}\label{am}
		A_m(\mu,\nu,w;z):=\sum_{k=0}^m\frac{(-1)^{-\mu-w-\frac{1}{2}}\Gamma\left(\mu+w+\frac{1}{2}+k\right)}{k!\Gamma\left(-\nu-\mu-k\right)\Gamma\left(\frac{1}{2}-\nu-\mu-w-k\right)}\left(\frac{z}{2}\right)^{-2k}.
	\end{align}
\end{theorem}
It is shown in \cite[Corollaries 2.6--2.9]{dkk} that Theorems \ref{befac} and \ref{extendedid} give in totality Ramanujan's formula \eqref{rameqn} as a corollary. 

Having general parameter $s$ also allows more flexibility. We can now let $s\to2m, m\geq0,$ in Theorem \ref{befac} and obtain \cite[Theorem 2.11]{dkk}:
\begin{theorem}\label{a=2mcase}
	 Let $\mathrm{Shi}(z)$ and $\mathrm{Chi}(z)$ be the hyperbolic sine and cosine integrals defined by \cite[p.~150, Equation (6.2.15), (6.2.16)]{olver-2010a}
	\begin{align}\label{shichi}
		\mathrm{Shi}(z):=\int_0^z\frac{\sinh(t)}{t}\ dt,\hspace{3mm}
		\mathrm{Chi}(z):=\gamma+\log(z)+\int_0^z\frac{\cosh(t)-1}{t}\ dt,
	\end{align}
where $\g$ is Euler's constant. Let $m\in\mathbb{N}$. Then for $\mathrm{Re}(y)>0$,
	\begin{align}\label{a=2midentity} 
		&\sum_{n=1}^\infty \sigma_{2m}(n)e^{-ny}-\frac{(2m)!}{y^{2m+1}}\zeta(2m+1)+\frac{B_{2m}}{2my}\nonumber\\
		&=(-1)^m\frac{2}{\pi}\left(\frac{2\pi}{y}\right)^{2m+1}\sum_{n=1}^\infty\sigma_{2m}(n)\Bigg\{\sinh\left(\frac{4\pi^2n}{y}\right)\mathrm{Shi}\left(\frac{4\pi^2n}{y}\right)\nonumber\\
		&\quad-\cosh\left(\frac{4\pi^2n}{y}\right)\mathrm{Chi}\left(\frac{4\pi^2n}{y}\right)+\sum_{j=1}^m(2j-1)!\left(\frac{4\pi^2n}{y}\right)^{-2j}\Bigg\}.
	\end{align}
\end{theorem}
Note that the series on the right-hand side of \eqref{a=2midentity} is a natural analogue of
\begin{equation*}
	\sum_{n=1}^\infty\sigma_{2m+1}(n)\Bigg\{\sinh\left(\frac{4\pi^2n}{y}\right)-\cosh\left(\frac{4\pi^2n}{y}\right)\Bigg\}=-\sum_{n=1}^\infty\sigma_{2m+1}(n)e^{-4\pi^2n/y},
\end{equation*}
which is one of the series appearing in the modular transformation satisfied by $\sum_{n=1}^{\infty}\sigma_{2m+1}(n)e^{-ny}$, that is, in \eqref{es}.

The case $m=0$ of Theorem \ref{a=2mcase} was obtained in an equivalent form by Wigert \cite[p.~203, Equation (A)]{wig0} (see also \cite[Equation (2.19)]{dkk}) who termed it \emph{`la formule importante'}.

Also, if we let $s\to-2m, m\in\mathbb{N}$, in Theorem \ref{extendedid}, we are led to \cite[Theorem 2.12]{dkk}:
\begin{theorem}\label{trans-2m}
Let $m\in\mathbb{N}$. If $\a$ and $\b$ are complex numbers such that $\textup{Re}(\a)>0$, $\textup{Re}(\b)>0$, and $\a\b=\pi^2$, then
	\begin{align}\label{equiforma=-2mano}
		&\a^{-\left(m-\frac{1}{2}\right)}\left\{\frac{1}{2}\zeta(2m)+\sum_{n=1}^\infty\frac{n^{-2m}}{e^{2n\a}-1}\right\}-\sum_{k=0}^{m-1}\frac{2^{2k-1}B_{2k}}{(2k)!}\zeta(2m-2k+1)\a^{2k-m-\frac{1}{2}}\nonumber\\
		&=(-1)^{m+1}\b^{-\left(m-\frac{1}{2}\right)}\left\{\frac{\g}{\pi}\zeta(2m)+\frac{1}{2\pi}\sum_{n=1}^\infty n^{-2m}\left(\psi\left(\frac{in\b}{\pi}\right)+\psi\left(-\frac{in\beta}{\pi}\right)\right)\right\},
	\end{align}
where $\psi(z)$ denotes the logarithmic derivative of $\Gamma(z)$.
\end{theorem}
 This result can be conceived as a \emph{companion} of Ramanujan's formula for $\zeta(2m+1)$. It was for the first time\footnote{Around the same time, Dorigoni and Kleinschmidt \cite{dorigoni} considered the case when $a$ is \emph{negative} even integer (see \cite[Equation (2.43)]{dorigoni}) using the concept of transseries \cite{transseries}.} in \cite{dkk} that the non-modular but explicit transformations given in Theorem \ref{befac} and \ref{extendedid} were obtained for $\sum_{n=1}^{\infty}\sigma_{2m}(n)e^{-ny}, m\in\mathbb{Z}\backslash\{0\}$. See \cite[p.~7]{dkk} for more discussion on these results. One of their applications is that one can easily obtain the asymptotic expansion of $\sum_{n=1}^{\infty}\sigma_{2m}(n)e^{-ny}$ using them. For example, for  $m\in\mathbb{N}$, we have \cite[Corollary 1.5]{dk03} as $y\to0$ in $|\arg(y)|<\pi/2$, 
 \begin{align}\label{asym for general m eqn}
 	\sum_{n=1}^\infty \sigma_{2m}(n)e^{-ny}&=\frac{(2m)!}{y^{2m+1}}\zeta(2m+1)-\frac{B_{2m}}{2my}-\frac{2(-1)^m}{\pi(2\pi)^{2m-1}}\sum_{j=1}^{r+1}\frac{\Gamma(2m+2j)\zeta(2m+2j)\zeta(2j)}{(2\pi)^{4j}}y^{2j-1}+O\left(y^{2r+3}\right).
 \end{align}
Ramanujan's first letter to Hardy contains the case $m=1$ of this asymptotic expansion and is proved by Watson using the Abel-Plana summation formula. This asymptotic expansion very easily leads to Wright's asymptotic estimate for the generating function of the number of plane partitions of a positive integer. For more details, see \cite[Corollary 1.6]{dk03} and the discussion preceding it. 
 
Another application of Theorem \ref{befac} is that the parameter $s$ being complex permits differentiation of both sides of the transformation. Before we present the identity resulting through this process, let us define $\psi_k(z)$ for $z\in\mathbb{C}\backslash\{a:a\in\mathbb{R}^{-}\cup\{0\}\}$ by
 \begin{equation}\label{psikx}
 	\psi_k(z):=-\gamma_k-\frac{\log^{k}(z)}{z}-\sum_{n=1}^{\infty}\left(\frac{\log^{k}(n+z)}{n+z}-\frac{\log^{k}(n)}{n}\right),
 \end{equation}
 where $\gamma_k$ is the generalized Stieltjes constant defined by \cite{berndthurwitzzeta}
 \begin{equation}\label{scz}
 	\gamma_k:=\lim_{n\to\infty}\left(\sum_{j=0}^{n}\frac{\log^{k}(j+1)}{j+1}-\frac{\log^{k+1}(n+1)}{k+1}\right).
 \end{equation}
  The $\psi_k(z)$ are the logarithmic derivatives of the generalized gamma functions of Dilcher \cite{dilcher}, and are intimately connected with the Laurent series coefficients of the Hurwitz zeta function $\zeta(s, z)$ around $s=1$. See \cite{dgs1} not only for a discussion on this but also for the literature survey.
 
 Differentiating \eqref{maineqn} with respect to $s$ and then letting $s\to0$ gives the following identity of Banerjee, Gupta and the current author \cite[Theorem 1.1]{bdg_log}:
 	\begin{theorem}\label{loglamb}
 	Let $\psi_1(z)$ be defined by
 	\eqref{psikx}. Then for $\textup{Re}(y)>0$,
 		\begin{align}\label{translog}
 			\sum_{n=1}^{\infty}\frac{\log(n) }{e^{ny}-1}&=-\frac{1}{4}\log(2\pi)+\frac{1}{2y}\log^{2}(y)-\frac{\gamma^2}{2y}+\frac{\pi^2}{12y}\nonumber\\
 			&\quad-\frac{2}{y}(\gamma+\log(y))\sum_{n=1}^{\infty}\left\{\log\left(\frac{2\pi n}{y}\right)-\frac{1}{2}\left(\psi\left(\frac{2\pi in}{y}\right)+\psi\left(-\frac{2\pi in}{y}\right)\right)\right\}\nonumber\\
 			&\quad+\frac{1}{y}\sum_{n=1}^{\infty}\left\{\psi_1\left(\frac{2\pi in}{y}\right)+\psi_1\left(-\frac{2\pi in}{y}\right)-\frac{1}{2}\left(\log^{2}\left(\frac{2\pi in}{y}\right)+\log^{2}\left(-\frac{2\pi in}{y}\right)\right)+\frac{y}{4n}\right\}.
 		\end{align}
 	\end{theorem}
  Prior to the discovery of \eqref{translog}, no information about $\displaystyle\sum_{n=1}^{\infty}\frac{\log(n) }{e^{ny}-1}$ was available, not even its asymptotic expansion as $y\to0$. The above exact formula for the series readily gives its asymptotic expansion as $y\to 0$ as can be seen from \cite[Theorem 1.2]{bdg_log}. The latter is shown to have an application in deriving the asymptotic expansion to all orders of the smoothly weighted moment of $\zeta\left(\frac{1}{2}-it\right)\zeta'\left(\frac{1}{2}+it\right)$ on the critical line, that is, of $ \displaystyle\int_{0}^{\infty}\zeta\left(\frac{1}{2}-it\right)\zeta'\left(\frac{1}{2}+it\right)e^{-\delta t}\, dt$ as $\delta\to0$; see Theorem 1.3 of \cite{bdg_log}.

We note that for Re$(s)>2$, the series in \eqref{rames} was also considered by Ramanujan himself  \cite[p.~269]{ramnote} who obtained the following beautiful transformation for it.
\begin{theorem}
Let $\alpha$ and $\beta$ be two positive real numbers such that $\alpha\beta=4\pi^2$. Then for $\mathrm{Re}(s)>2$, we have 
\begin{align}\label{ramanujan gen eqn}
	\alpha^{s/2}\left\{\frac{\Gamma(s)\zeta(s)}{(2\pi)^s}+\cos\left(\frac{\pi s}{2}\right)\sum_{n=1}^\infty\frac{n^{s-1}}{e^{n\alpha}-1}\right\}&=\beta^{s/2}\left\{\cos\left(\frac{\pi s}{2}\right)\frac{\Gamma(s)\zeta(s)}{(2\pi)^s}+\sum_{n=1}^\infty\frac{n^{s-1}}{e^{n\beta}-1}\right.\nonumber\\
	&\left.\quad-\sin\left(\frac{\pi s}{2}\right)\mathrm{PV}\int_0^\infty\frac{x^{s-1}}{e^{2\pi x}-1}\cot\left(\frac{1}{2}\beta x\right)dx\right\},
\end{align}
where $\mathrm{PV}$ denotes the Cauchy principal value integral. 
\end{theorem}
It is clear that letting $s=2\ell$ and respectively replacing $\alpha$ and $\beta$ by $2\alpha$ and $2\beta$ in \eqref{ramanujan gen eqn} gives \eqref{es}, and hence \eqref{ramanujan gen eqn} is also a \emph{continuous version} of \eqref{es}. Berndt proved the above formula in \cite[p.~416]{bcbramfifthnote} using the Abel-Plana summation formula. 

Two different generalizations of this result have been obtained recently. The first one is given by Berndt, Gupta, Zaharescu and the current author \cite[Theorem 14]{bdgz} and is in the setting of Koshliakov zeta functions which are shortly discussed below. The second generalization was obtained by Kumar and the current author \cite[Theorem 1.1]{dk03} and is given next.
\begin{theorem}\label{ram with a}
	Let $\mathrm{Re}(\alpha),\mathrm{Re}(\beta)>0$ such that $\alpha\beta=4\pi^2$. Let $0\leq a<1$. Then, for $\mathrm{Re}(s)>2$, the following transformation holds:
	\begin{align}\label{ram with a eqn}
		&\alpha^{s/2}\left\{\frac{\Gamma(s)\zeta(s)}{(2\pi)^s}+\frac{1}{2}\sum_{n=1}^\infty n^{s-1}\left(\frac{e^{\pi is/2}}{e^{n\alpha-2\pi ia}-1}+\frac{e^{-\pi is/2}}{e^{n\alpha+2\pi ia}-1}\right)\right\}\nonumber\\
		&=\beta^{s/2}\left\{\frac{\Gamma(s)}{(2\pi)^{s}}\sum_{k=1}^\infty\frac{\cos\left(\frac{\pi s}{2}+2\pi ak\right)}{k^s}+\sum_{n=1}^\infty\frac{(n-a)^{s-1}}{e^{(n-a)\beta}-1}\right.\nonumber\\
		&\quad\left.-\frac{1}{2i}\mathrm{PV}\int_0^\infty x^{s-1}
		\left(\frac{e^{\pi is/2}}{e^{2\pi x-2\pi ia}-1}-\frac{e^{-\pi is/2}}{e^{2\pi x+2\pi ia}-1}\right)\cot\left(\frac{1}{2}\beta x\right)dx\right\}.
	\end{align}
	\end{theorem}
The above result involves the generalized Lambert series 
\begin{equation*}
	\sum_{n=1}^{\infty}\frac{(n-a)^{s-1}}{e^{(n-a)z}-1}\hspace{8mm}(s\in\mathbb{C}, \textup{Re}(z)>0, 0\leq a<1),
\end{equation*}
which has not been studied before except for some special values such as $a=1/2$ or $1/4$. Three special cases of \eqref{ram with a eqn} are derived in \cite[Corollaries 1.2-1.4]{dk03}.\\

\subsection{A generalization in the setting of Koshliakov zeta functions}
Nikolai Sergeevich Koshliakov \cite{koshliakov3} wrote a beautiful manuscript in 1949 that lay dormant in the mathematical community for over 70 years. In this manuscript, he developed the theory of generalized zeta functions and the functions associated with them. The theory has its genesis in a problem on heat conduction resulting from Physics (see \cite[Section 2]{dg}).  This manuscript was studied in detail for the first time by Gupta and the current author in \cite{dg} who also built the theory further by obtaining two new modular equations, one of which is a new generalization of \eqref{rameqn}. It concerns one of the two \emph{Koshliakov zeta functions} defined below.

Let $p>0$.  
Then the first Koshliakov zeta function $\zeta_p(s)$ is defined by \cite[p.~6]{koshliakov3}
\begin{align}\label{zps}
	\zeta_p(s):=\sum_{j=1}^{\infty}\frac{p^2+\lambda^2_{j}}{p\left(p+\frac{1}{\pi}\right)+\lambda^2_j}\cdot \frac{1}{\lambda^s_j},\qquad
	\Re(s)>1,
\end{align}
where $\lambda_j$ runs over the roots of the transcendental equation $p \sin(\pi \lambda)+\lambda \cos(\pi \lambda)=0.$ Observe that $\lim_{p\to\infty}\zeta_p(s)=\zeta(s)$.

Let $\sigma(z):=\frac{p+z}{p-z}$ and $\sigma_p(z):=\sum_{j=1}^{\infty}\frac{p^2+\lambda_j^2}{p\left(p+\frac{1}{\pi}\right)+\lambda_j^2}e^{-\l_j z}$. For $k\in\mathbb{N}$, Koshliakov's generalized Bernoulli numbers\footnote{We note that in Koshliakov's notation, $B_{2k}^{(p)}$ would be denoted by $(-1)^{k+1}B_k^{(p)}$. We have followed the contemporary notation for Bernoulli numbers. It is easy to see that $\lim_{p\to\infty}B_{2k}^{(p)}=B_{2k}$.} are defined by \cite[p.~46, Chapter 2, Equation~(45)]{koshliakov3}
\begin{align}\label{genber}
	B_{2k}^{(p)} :=(-1)^{k+1}4k\int_{0}^{\infty}x^{2k-1}\sigma_p(2\pi x)dx,\hspace{1mm} B_{0}^{(p)}:=\frac{1}{1+\frac{1}{\pi p}}.
\end{align}
Then for $m \in \mathbb{Z},~m\neq 0,$ and $\a \b =\pi^2$, the generalization of Ramanujan's formula \eqref{rameqn} derived in \cite[Theorem 4.1]{dg} is
\begin{align}\label{rt1eqn}
	&\a^{-m}\left\{\frac{1}{2}\zeta_{p}(2m+1)+\sum_{j=1}^{\infty}\frac{p^2+\l_j^2}{p\left(p+\frac{1}{\pi}\right)+\l_j^2}\cdot\frac{\l_{j}^{-2m-1}}{\sigma\left(\frac{\l_j \a}{\pi} \right)e^{2\a \l_j}-1}\right\}\nonumber\\
	&=(-\b)^{-m}\left\{\frac{1}{2}\zeta_{p}(2m+1)+\sum_{j=1}^{\infty}\frac{p^2+\l_j^2}{p\left(p+\frac{1}{\pi}\right)+\l_j^2}\cdot\frac{\l_{j}^{-2m-1}}{\sigma\left(\frac{\l_j \b}{\pi} \right)e^{2\b \l_j}-1}\right\}\nonumber\\
	&\quad -2^{2m}\sum_{j=0}^{m+1}\frac{(-1)^jB_{2j}^{(p)}B_{2m-2j+2}^{(p)}}{(2j)!(2m-2j+2)!}\a^{m-j+1}\b^j.
\end{align}
It is clear that letting $p\to\infty$ yields \eqref{rameqn}. Also, if we let $p\to0$, we get an analogue of \eqref{rameqn}; see the papers of Malurkar \cite{malurkar}, Berndt \cite{berndtcrelle}, Gupta and the current author \cite[Corollary 4.2]{dg} and Chourasiya, Jamal and Maji \cite{cjm}.

\subsection{Generalizing the classical theory of Eisenstein series}

In 2001 \cite{ktyhr}, Kanemitsu, Tanigawa and Yashimoto studied the generalized Lambert series $\sum\limits_{n=1}^{\infty}\frac{n^{N-2h}}{e^{n^{N}x}-1}$ for $h, N\in\mathbb{N}$ such that $h\leq N/2$ and obtained a transformation for it \cite[Theorem 1]{ktyhr}. A slight generalization of this series was first considered by Ramanujan! See page 332 of Ramanujan's Lost Notebook \cite{lnb}. However, Ramanujan did not give any result for this series. A possibility of some pages in the Lost Notebook being lost cannot be discarded altogether.

In \cite[p.~385--386]{RLNII}, Andrews and Berndt discuss a bit on placing the series in a framework generalizing the classical theory of Eisenstein series. Samplings of such a theory have recently been conceived in the papers \cite{dixitmaji1}, \cite{dgkm} and \cite{BDG}. In \cite[Theorem 1.2]{dixitmaji1}, Maji and the current author give a generalization of Ramanujan's formula \eqref{rameqn} containing the aforementioned generalized Lambert series considered by Ramanujan. This formula is given next. The nice thing about it is that it gives a relation between two \emph{different} odd zeta values of the form $\zeta(2m+1)$ and $\zeta(2Nm+1)$, where $N$ is an odd natural number, by means of these Lambert series.
\begin{theorem}\label{dm1}
Let $N$ be an odd positive integer and $\a,\b>0$ such that $\a\b^{N}=\pi^{N+1}$. Then for $m\in\mathbb{Z}, m\neq 0$,
\begin{align}\label{zetageneqn}
	&\a^{-\frac{2Nm}{N+1}}\left(\frac{1}{2}\zeta(2Nm+1)+\sum_{n=1}^{\infty}\frac{n^{-2Nm-1}}{\textup{exp}\left((2n)^{N}\a\right)-1}\right)\nonumber\\
	&=\left(-\b^{\frac{2N}{N+1}}\right)^{-m}\frac{2^{2m(N-1)}}{N}\Bigg(\frac{1}{2}\zeta(2m+1)+(-1)^{\frac{N+3}{2}}\sum_{j=\frac{-(N-1)}{2}}^{\frac{N-1}{2}}(-1)^{j}\sum_{n=1}^{\infty}\frac{n^{-2m-1}}{\textup{exp}\left((2n)^{\frac{1}{N}}\b e^{\frac{i\pi j}{N}}\right)-1}\Bigg)\nonumber\\
	&\quad+(-1)^{m+\frac{N+3}{2}}2^{2Nm}\sum_{j=0}^{\left\lfloor\frac{N+1}{2N}+m\right\rfloor}\frac{(-1)^jB_{2j}B_{N+1+2N(m-j)}}{(2j)!(N+1+2N(m-j))!}\a^{\frac{2j}{N+1}}\b^{N+\frac{2N^2(m-j)}{N+1}}.
\end{align}
\end{theorem}
It is straightforward to see that letting $N=1$ in \eqref{zetageneqn} gives \eqref{rameqn}. A one-parameter generalization of a transformation equivalent to that of the Dedekind eta function, that is, of \eqref{etatransequiv}, is also derived in \cite[Theorem 1.3]{dixitmaji1}. The title of this subsection is in light of the fact that \eqref{rameqn} gives the classical modular transformations of the Eisenstein series. A character analogue of \eqref{zetageneqn} has been recently derived by Gupta, Jamal, Karak and Maji  \cite[Theorem 2.2]{gjkm}.

Two different directions were recently undertaken to further generalize \eqref{zetageneqn}. We begin with the one in \cite{dgkm} where Gupta, Kumar, Maji and the current author studied the more general Lambert series considered by Kanemitsu, Tanigawa and Yoshimoto in \cite{ktyacta}, namely, \begin{equation}\label{gls}
	\sum_{n=1}^{\infty}n^{N-2h}\frac{\textup{exp}(-an^{N}x)}{1-\textup{exp}(-n^{N}x)},
\end{equation}
and obtained a two-parameter generalization of \eqref{rameqn} containing this series. Their result is as follows.
\begin{theorem}\label{ggram}
	Let $0<a\leq 1$, let $N$ be an odd positive integer and $\a,\b>0$ such that $\a\b^{N}=\pi^{N+1}$. Then for any positive integer $m$,
	\begin{align}\label{zetageneqna}
		&\a^{-\frac{2Nm}{N+1}}\bigg(\left(a-\frac{1}{2}\right)\zeta(2Nm+1)+\sum_{j=1}^{m-1}\frac{B_{2j+1}(a)}{(2j+1)!}\zeta(2Nm+1-2jN)(2^N\a)^{2j}\nonumber\\
		&\qquad\qquad+\sum_{n=1}^{\infty}\frac{n^{-2Nm-1}\textup{exp}\left(-a(2n)^{N}\a\right)}{1-\textup{exp}\left(-(2n)^{N}\a\right)}\bigg)\nonumber\\
		&=\left(-\b^{\frac{2N}{N+1}}\right)^{-m}\frac{2^{2m(N-1)}}{N}\bigg[\frac{(-1)^{m+1}(2\pi)^{2m}B_{2m+1}(a)N \g}{(2m+1)!}+\frac{1}{2}\sum_{n=1}^{\infty}\frac{\cos(2\pi na)}{n^{2m+1}}\nonumber\\
		&\quad+(-1)^{\frac{N+3}{2}}\sum_{j=\frac{-(N-1)}{2}}^{\frac{N-1}{2}}(-1)^{j}\bigg\{\sum_{n=1}^{\infty}\frac{n^{-2m-1}\cos(2\pi na)}{\textup{exp}\left((2n)^{\frac{1}{N}}\b e^{\frac{i\pi j}{N}}\right)-1}\nonumber\\
		&\quad+\frac{(-1)^{j+\frac{N+3}{2}}}{2\pi}\sum_{n=1}^{\infty}\frac{\sin(2\pi na)}{n^{2m+1}}\left(\psi\left(\tfrac{i\beta}{2\pi}(2n)^{\frac{1}{N}} e^{\frac{i\pi j}{N}}\right)+\psi\left(\tfrac{-i\beta}{2\pi}(2n)^{\frac{1}{N}} e^{\frac{i\pi j}{N}}\right)\right)\bigg\}\bigg]\nonumber\\
		&\quad+(-1)^{m+\frac{N+3}{2}}2^{2Nm}\sum_{j=0}^{\left\lfloor\frac{N+1}{2N}+m\right\rfloor}\frac{(-1)^jB_{2j}(a)B_{N+1+2N(m-j)}}{(2j)!(N+1+2N(m-j))!}\a^{\frac{2j}{N+1}}\b^{N+\frac{2N^2(m-j)}{N+1}}.
	\end{align}
\end{theorem}
It is not difficult to see that letting $a=1$ in Theorem \ref{ggram} gives Theorem \ref{dm1} for $m\in\mathbb{N}$. Observe that through \eqref{zetageneqna}, we can provide a relation between \emph{any} number of odd zeta values $\zeta(2N+1)$, $\zeta(4N+1)$, $\zeta(6N+1), \cdots, \zeta(2Nm+1)$, where  $N$ is an odd positive integer  and $m\in\mathbb{N}$, in terms of generalized Lambert series and generalized higher Herglotz functions (see \eqref{vlaszag0} below for the definition of higher Herglotz functions). Numerous other corollaries resulting from Theorem \ref{ggram} are also given in \cite{dgkm}. For example, in \cite[Corollary 2.5]{dgkm}, a transformation linking $\zeta(3)/\pi^3$, $\zeta(5)/\pi^5$, $\zeta(7)/\pi^7$, $\zeta(9)/\pi^9$ and $\zeta(11)/\pi^{11}$ is obtained.

A generalization of \eqref{extendedid}, with an extra parameter $N$ was recently derived by Banerjee, Gupta and the current author \cite[Theorem 2.3]{BDG}.

\begin{theorem}\label{Analytic continuation}
	Let $m \in \mathbb{N} \cup \lbrace 0 \rbrace$ and \textup{Re}$(y)>0$. For $s\in\mathbb{C}$ and $N\in\mathbb{N}$, let
	\begin{equation}\label{defaf}
		\sigma_s^{(N)}(n) := \sum_{d^N\mid n} d^s,\hspace{6mm}S_{s}^{(N)}(n):= \sum\limits_{d_{1}^N d_2=n} d_{2}^{\frac{1+s}{N}-1}.
	\end{equation} 
For \textup{Re}$(s)>-(2m+2)N-1$, the following identity holds:
\begin{align}\label{Eqn:Analytic continuation}
	&\sum_{n=1}^\infty \sigma_s^{(N)}(n) e^{-ny}  + \frac{\zeta(-s)}{2} - \frac{\zeta(N-s)}{y} - \frac{1}{N} \frac{\Gamma \left(\frac{1+s}{N}\right) \zeta \left(\frac{1+s}{N}\right)}{y^{\frac{1+s}{N}}} \nonumber\\
	&=\frac{y}{2\pi^{2}} \sum_{k=0}^m \left(-\frac{y^2}{4\pi^2}\right)^k   \zeta(-2kN-N-s)\zeta(2k+2)+\frac{2(2\pi)^{\frac{1}{N}-\frac{1}{2}}N^{\frac{s-1}{2}}}{y^{\frac{1}{N}+\frac{s}{2N}}}\nonumber\\
	&\quad\times\sum_{n=1}^{\infty}\frac{S_{s}^{(N)}(n)}{n^{\frac{s}{2N}}}\bigg[{}_\frac{1}{2}K_{\frac{s}{2N}}^{(N)}(n)\left(\frac{4\pi^{N+1}n}{yN^{N}}, 0\right)-\frac{2^{\frac{1}{2}+\frac{s+1}{N}}\pi^{\frac{(1-N)a}{2N}-N}}{\left(\frac{4\pi^{N+1}n}{yN^{N}}\right)^{1+\frac{1}{N}+\frac{s}{2N}}}\frac{\sin\left(\frac{\pi}{2}(N-s)\right)}{2^{N-1}}C_{m, N}\left(\frac{1}{2}, \frac{s}{2N}, 0, \frac{4\pi^{N+1}n}{yN^{N}}\right)\bigg],
\end{align}
	where
		\begin{equation}\label{muknug}
		{}_{\mu}K_{\nu}^{(N)}(z, w) := 2^{\mu + \frac{2}{N}-1}\pi^{(1-N)\nu}z^{w+\nu-\frac{2}{N}} G_{1, \, \, 2N+1}^{N+1, \, \, 1}\bigg( \begin{matrix}
			1+\frac{1}{2N}-\mu-\nu-w \\
			\frac{1}{2}+\frac{1}{2N}-\nu, \langle \frac{i}{N} \rangle_{i=1}^N; 1+\frac{1}{2N}-w, \langle 1+\frac{3}{2N} - \frac{i}{N} \rangle_{i=2}^N 
		\end{matrix} \bigg | \frac{z^2}{4} \bigg),
	\end{equation}
with $G_{p,q}^{\,m,n} \!\left(  \,\begin{matrix} a_1,\cdots , a_p \\ b_1, \cdots b_m; b_{m+1}, \cdots, b_q \end{matrix} \; \Big| X   \right) $ being the Meijer G-function \cite[p.~415, Definition 16.17]{olver-2010a}, and 
	\begin{align}\label{cmn}
		C_{m, N}(\mu, \nu, w; z) &:= \sum_{k=0}^m \frac{(-1)^{k(N+1)+N}}{k!}  \Gamma \left(\tfrac{1}{2}+\mu+w +k\right) \Gamma\left(1+\mu+\nu +k\right)\prod_{i=1}^{2N-1}\Gamma\left( \tfrac{i}{2N} +\mu+\nu+w+k\right)\left(\frac{z}{2} \right)^{-2k}.
	\end{align}
\end{theorem}
Letting $N=1$ in the above theorem gives Theorem \ref{extendedid}. There are several important  corollaries of the above result. For example, when $s=-2Nm-N,$ where $N$ is an odd positive integer, we get a new generalization of \eqref{rameqn}.

 We note in passing that Gupta and Maji \cite[Corollary 3.5]{guptamaji} have obtained another generalization of \eqref{rameqn}. In fact, they further generalize their result \cite[Theorem 3.1]{guptamaji}. See also the last paragraph of \cite[Section 11]{BDG}.

\subsection{Ramanujan's formula through period polynomials and its generalization}
As mentioned in antepenultimate paragraph of the introduction, Ramanujan's formula can be rephrased in terms of a representation for the period polynomials of the Eisenstein series on the full modular group. If we consider analogous period polynomials for the Eisenstein series of higher level, that would then lead us to Ramanujan-type formulas. This was done by Berndt and Straub in \cite[Theorem 6.1]{berndtstraubmathz} for a generalized Eisenstein series associated with primitive Dirichlet characters $\chi$ and $\psi$ modulo $L$ and $M$ respectively. In particular, when $\chi\equiv1$ and $\psi=\chi_{-4}$, the non-principal Dirichlet character modulo 4, defined by $\chi_{-4}(n)=0$ for even $n$ and $\chi_{-4}(n)=(-1)^{(n-1)/2}$ for odd $n$, then the special case of their result is an analogue of \eqref{rameqn} given by Ramanujan himself \cite[p.~277, Entry 21(iii)]{bcbramsecnote}:
\begin{theorem}
	Let $\alpha, \beta>0$ with $\alpha\beta=\pi^2$. Let $L(\psi, k):=\sum_{n=1}^{\infty}\psi(n)n^{-k}$ for $k>1$. For $m\in\mathbb{N}$,
	\begin{align}\label{entry21iii}
		\alpha^{-m+1/2}\left\{\frac{1}{2}L(\chi_{-4}, 2m)+\sum_{n=1}^{\infty}\frac{\chi_{-4}(n)n^{-2m}}{e^{n\alpha}-1}\right\}&=\frac{(-1)^{m}\beta^{-m+1/2}}{2^{2m+1}}\sum_{n=1}^{\infty}\frac{\textup{sech}(n\beta)}{n^{2m}}\nonumber\\
		&\quad+\frac{1}{4}\sum_{n=0}^{m}\frac{(-1)^n}{2^{2n}}\frac{E_{2n}}{(2n)!}\frac{B_{2m-2n}}{(2m-2n)!}\alpha^{m-n}\beta^{n+1/2},
	\end{align}  
	where $E_j$ denotes the $j^{\textup{th}}$ Euler number \cite[p.~15]{temme}.
\end{theorem}

We note in passing that a special case of \eqref{ram with a eqn} derived in \cite[Corollary 1.4]{dk03} gives an integral representation for a cousin of the Lambert series occurring in \eqref{entry21iii}, namely, $\sum\limits_{n=1}^{\infty}\displaystyle\frac{\chi_{-4}(n)n^{2m-1}}{e^{n\beta}-1}$, where $m\in\mathbb{N}, m>1$. The latter does not fall under the purview of the setting in \cite{berndtstraubmathz}.
%For a fixed $j\in\mathbb{N}$, let $P_j(X_1, \cdots, X_r)$ denote a polynomial of $r$ variables such that $P_j(n)>0$ for all $n\in\mathbb{N}^{r}$. Then the  series
%\begin{equation}
%	\sum_{n\in\mathbb{N}^{r}}P_j(n)^{-s}=\sum_{n_1=1}^{\infty}\cdots\sum_{n_r=1}^{\infty}P_j(n_1, \cdots, n_r)^{-s}
%\end{equation}
%converges absolutely for Re$(s)>\sigma_j>0$. Let $\beta=(\beta_1, \cdots, \beta_r)$ denote an $r$-tuple of non-negative integers. In \cite{eie-chen}, Eie and Chen consider the zeta function $Z(P_j, \beta,s)$ associated with $P_j$ defined by
%\begin{align}
%	Z(P_j, \beta, s):=\sum_{n\in\mathbb{N}^{r}}n^{\beta}P_j(n)^{-s}
%	=\sum_{n_1=1}^{\infty}\cdots\sum_{n_r=1}^{\infty}n_1^{\beta_1}\cdots n_{r}^{\beta_r}P_j(n_1, \cdots, n_r)^{-s},
%\end{align}
%where Re$(s)>\sigma_j+|\beta|$ with $|\beta|=\beta_1+\cdots+\beta_r$. They prove that $Z(P_j, \beta; s)$ as well as $Z(\prod_{j=1}^n P_j, \beta; s)$ have meromorphic continuations in the whole complex plane. Furthermore, they prove that
%\begin{equation*}
%	Z\left(\prod_{j=1}^n P_j, \beta; 0\right)=\frac{1}{n}\sum_{j=1}^{n}Z(P_j, \beta; 0).
%\end{equation*}
%As an immediate application, they use the above relation to evaluate the special values of zeta functions associated with products of linear forms previously considered by Eie as well as by Shintani.
%
%As another application of their result, they obtain a beautiful generalization as well as an analogue of Ramanujan's formula \eqref{rameqn}; see Proposition 3 of \cite{eie-chen}.\

\subsection{Other generalizations}

Franke \cite{franke-rj} has extended \eqref{rameqn} in two different ways by considering generalized Dirichlet series having properties similar to the Dirichlet $L$-functions. Also see \cite{franke-rint}.
Recently, Banerjee, Gupta and Kumar \cite[Theorem 1.1]{bgk} have obtained a generalization of \eqref{rameqn} in the setting of the Dedekind zeta function over an arbitrary number field. 
Ramanujan's theorem has been massively generalized by Katsurada \cite{katsurada} and Lim \cite[Theorem 1.1]{lim1}.
Chavan, Chavan, Vignat and Wakhre \cite{ccvw} have shown that Ramanujan's formula is a special case of their more general result \cite[Theorem 2.2]{ccvw} on the convolution of generalized Dirichlet series \cite[Equation (2.3)]{ccvw} parametrized by a set of zeros and certain weights. Kongsiriwong \cite[Theorem 2.5]{kongsiriwong} has given yet another extension of \eqref{rameqn} using certain infinite series involving cotangent functions. An identity wherein the product of Riemann zeta functions in the finite sum on the right-hand side of \eqref{rameqn} is replaced by the corresponding one involving Hurwitz zeta functions has been derived by Chavan \cite{pchavan}. Ramanujan's formula follows as a special case of his identity.

\section{Analogues of Ramanujan's formula \eqref{rameqn}}

\subsection{A non-holomorphic counterpart of Ramanujan's formula}

O'Sullivan \cite[Theorem 1.3]{sullivan} has recently obtained a beautiful non-holomorphic analogue of \eqref{rameqn}:
\begin{theorem}\label{non-holomorphic-thm}
Let $z=x+iy\in\mathbb{H}$. Define $V_k(z)$ by 
\begin{align}\label{vkz}
	V_k(z):=\sum_{n=1}^{\infty}\sigma_{k-1}(n)e^{-2\pi in\bar{z}}\sum_{u=0}^{-k}\frac{(4\pi ny)^u}{u!}.
\end{align}
For all $k\in2\mathbb{Z}$,
\begin{align}\label{non-holomorphic}
	2\left(z^kV_k(z)-V_k(-1/z)\right)&=\frac{2\zeta(2-k)}{(2\pi i)^k}\left(\frac{y}{\pi}\right)^{1-k}\left(|z|^{2k-2}-z^k\right)\nonumber\\
	&\quad-(2\pi i)^{1-k}\sum_{\substack{u, v\in\mathbb{Z}_{\geq0}\\u+v=1-k/2}}\frac{B_{2u}}{(2u)!}\frac{B_{2v}}{(2v)!}z^{1-2v}+\begin{cases}
		0,\hspace{32.5mm}\text{if}\hspace{1mm}k>0,\\
		\pi i/2+\overline{\log(z)}, \hspace{11mm}\text{if}\hspace{1mm}k=0,\\
		(1-z^k)\zeta(1-k),\hspace{6mm}\text{if}\hspace{1mm}k<0.
\end{cases}\end{align}	
\end{theorem}

\subsection{A Ramanujan-type formula involving the higher Herglotz functions}
For $k\in\mathbb{N},k>1,$ and $x\in\mathbb{C}\backslash(-\infty,0]$, let $F_k(x)$ be the higher Herglotz function defined by  \begin{equation}\label{vlaszag0}
	F_k(x):=\sum_{n=1}^{\infty}\frac{\psi(nx)}{n^{k}}.
\end{equation}
This function has played an important role in the work of Vlasenko and Zagier \cite{vz} on deriving higher Kronecker ``limit" formulas for real quadratic fields. It was first studied by Maier \cite[p.~114]{maier}. 

In \cite[Corollary 3.4]{dgk}, Gupta, Kumar and the current author obtained an analogue of Ramanujan's formula involving $F_k(x)$ given below\footnote{They obtain, in fact, a result more general than \eqref{trans2m+1}; see \cite[Theorem 3.2]{dgk}.}.
\begin{theorem}\label{trans2m+1}
	Let $\alpha$ and $\beta$ be two complex numbers with $\textup{Re}(\a)>0, \textup{Re}(\b)>0$ and $\alpha\beta=\pi^2$. Then for $m\in\mathbb{N}$,
	\begin{align}\label{trans2m+1eqn}
		&\alpha^{-m}\left\{2\gamma\zeta(2m+1)+F_{2m+1}\left(\frac{i\alpha}{ \pi}\right)+F_{2m+1}\left(-\frac{i\alpha}{\pi}\right)\right\}\nonumber\\
		&=-(-\beta)^{-m}\left\{2\gamma\zeta(2m+1)+F_{2m+1}\left(\frac{i\beta}{ \pi}\right)+F_{2m+1}\left(-\frac{i\beta}{ \pi}\right)\right\}\nonumber\\
		&\qquad-2\sum_{j=1}^{m-1}(-1)^j\zeta(1-2j+2m)\zeta(2j+1)\alpha^{j-m}\beta^{-j}.
	\end{align}
\end{theorem}
In particular, letting $m=1$ in Theorem \ref{trans2m+1} gives the beautiful modular relation \cite[Corollary 3.5]{dgk}:
\begin{align}\label{modulareqn}
	\frac{1}{\a}\left\{2\g\zeta(3)+F_{3}\left(\frac{i\alpha}{2 \pi}\right)+F_{3}\left(-\frac{i\alpha}{2 \pi}\right)\right\}=\frac{1}{\b}\left\{2\g\zeta(3)+F_{3}\left(\frac{i\beta}{2 \pi}\right)+F_{3}\left(-\frac{i\beta}{2 \pi}\right)\right\}.
\end{align}
Observe the three different combinations - two even zeta values, one even and one odd zeta value, and two odd zeta values - occurring in the finite sums on the right-hand sides of \eqref{rameqn}, \eqref{equiforma=-2mano} and \eqref{trans2m+1eqn} respectively upon using Euler's formula \eqref{zetaevenint} to write $B_{2k}$ in terms of $\zeta(2k)$. 

\subsection{A Ramanujan-type  formula for $\zeta^2(2m+1)$}
Recently, Gupta and the current author \cite[Theorem 2.1]{dgsquare} obtained a Ramanujan-type formula for $\zeta^2(2m+1)$ stated below.
\begin{theorem}\label{zetasquared}
	Let $\epsilon=e^{i\pi/4}$ and $\overline{\epsilon}=e^{-i\pi/4}$. For $\rho>0$ and $x>0$, define $\Omega_{\rho}(x)$ by
	\begin{equation*}
		\Omega_{\rho}(x):=2\sum_{j=1}^{\infty} d(j)\left(K_{0}(4\rho\epsilon\sqrt{jx})+K_{0}(4\rho\overline{\epsilon}\sqrt{jx})\right),
	\end{equation*}
where $K_0(x)$ denotes the modified Bessel function of the second kind of order zero.
	Let $m$ be a non-zero integer. Define $\mathfrak{F}_{m}(\rho)$ by
	\begin{align}\label{zetasquaredeqn0}
		\mathfrak{F}_{m}(\rho):=(\rho^2)^{-m} \left \{ \zeta^{2}(2m+1)\left(\gamma + \log\left(\frac{\rho}{\pi}\right)- \frac{\zeta'(2m+1)}{\zeta(2m+1)}\right)+\sum_{n=1}^{\infty} \frac{d(n)\Omega_{\rho}(n)}{n^{2m+1}} \right \}.
	\end{align}
	Then for any $\alpha, \beta >0$ satisfying $\alpha \beta = \pi^2$,
	\begin{align}\label{zetasquaredeqn}
		\mathfrak{F}_{m}(\a)- (-1)^{-m}\mathfrak{F}_{m}(\b)=- \pi 2^{4m}\sum_{j=0}^{m+1}\frac{(-1)^{j}B^2_{2j}B^2_{2m+2-2j}}{(2j)!^2(2m+2-2j)!^2}\alpha^{2j}
		\beta^{2m+2-2j}.
	\end{align}
\end{theorem}
Two generalizations of Theorem \ref{zetasquared} are also derived in the same paper; see \cite[Theorems 2.2, 5.1]{dgsquare}. A common generalization of \eqref{rameqn} and Theorem \ref{zetasquared} has been given by Banerjee and Sahani \cite{banerjee-sahani} by obtaining a Ramanujan-type formula for $\zeta^{k}(2m+1)$ for $k\in\mathbb{N}$.\\

\subsection{A multidimensional analogue of an identity of Ramanujan}
Ramanujan has derived another transformation related to \eqref{entry21iii} (see \cite[p.~276, Entry 21(ii)]{bcbramsecnote}). As a special case, it gives 
\begin{align}\label{befmul}
	\sum_{n=1}^{\infty}\frac{\chi_{-4}(n)}{n}\textup{sech}\left(\frac{\pi n}{2}\right)=\frac{\pi}{8}.
\end{align}
A multidimensional analogue of \eqref{befmul} was recently obtained by Daniyarkhodzhaev and Korolev \cite[Theorem 1]{danbiyarkhodzhaev}. It states\footnote{The formula, as stated in \cite{danbiyarkhodzhaev} has a typo, namely, the occurrence of $\pi/2$ inside the argument of $\cosh$ should be replaced by $\pi$.} that for any $r\geq1$, 
\begin{equation*}
\sum_{n_1,\cdots, n_r=0}^{\infty}\frac{(-1)^{n_1+\cdots+n_r}}{\left(n_1+\frac{1}{2}\right)\cdots\left(n_r+\frac{1}{2}\right)\cosh\left(\pi\sqrt{\left(n_1+\frac{1}{2}\right)^2+\cdots+\left(n_r+\frac{1}{2}\right)^2}\right)}=\frac{1}{r+1}	\left(\frac{\pi}{2}\right)^r.
\end{equation*}

\subsection{Ramanujan polynomials and their cousins}
In \cite{msw}, Murty, Smyth and Wang defined \emph{Ramanujan polynomial} to be\footnote{Gun, Murty and Rath \cite{gmr} have also defined these polynomials although in their definition the power of $z$ is $2m+2-2k$ rather than $2k$. The reciprocal property in \eqref{rp} is seen to hold in both the definitions.}
\begin{align}\label{rp}
	R_{2m+1}(z):=\sum_{k=0}^{m+1}\frac{B_{2k}B_{2m+2-2k}}{(2k)!(2m+2-2k)!}z^{2k},
\end{align}
which is clearly constructed from the finite sum on the right-hand side of Ramanujan's formula \eqref{rameqn} upon replacing $k$ by $m+1-k$, and then letting $\alpha=-i\pi z$ and $\beta=i\pi/z$. Ramanujan polynomials have nice properties: they are reciprocal polynomials, that is, they satisfy the functional equation
\begin{equation*}
	R_{2m+1}(z)=z^{2m+2}R_{2m+1}\left(\frac{1}{z}\right),
\end{equation*}
they have real coefficients, and moreover, all of their non-real zeros lie on the unit circle etc. In the language of modular forms, these polynomials can be thought in terms of period polynomials of the Eisenstein series \cite[p.~4762]{conrey}. Other generalizations of the Ramanujan polynomials are considered in \cite{lalin-rogers}, \cite{lalin-smyth}, \cite[Equations (43), (47)]{berndtstraubmathz} and \cite{jmos}.

We now explain the importance of studying the zeros of Ramanujan polynomials as described in \cite{msw}. Letting $\alpha=-i\pi z$ and $\beta=i\pi/z$, where Im$(z)>0$, $F_{k}(z):=\sum_{n=1}^{\infty}\sigma_k(n)n^{-k}e^{2\pi inz}$ in \eqref{rameqn} with $m\in\mathbb{N}$, we get \cite{gross1}
\begin{align}
	F_{2m+1}(z)-z^{2m}F_{2m+1}\left(-\frac{1}{z}\right)	=\frac{1}{2}(z^{2m}-1)\zeta(2m+1)+\frac{1}{2z}(2\pi i)^{2m+1}R_{2m+1}(z).
\end{align}
Now Murty, Smyth and Wang \cite{msw} have proved that for each $m\geq 4$, there exists at least one algebraic number $\kappa$ (depending on $m$) with $|\kappa|=1, \kappa^{2m}\neq 1$, and lying in the upper half plane such that $R_{2m+1}(\kappa)=0$. Hence we get the following formula which expresses the odd zeta values $\zeta(2m+1), m\geq 4,$ in terms of two Eichler integrals:
\begin{align}
	\zeta(2m+1)=\frac{2}{\kappa^{2m}-1}\left(F_{2m+1}(\kappa)-\kappa^{2m}F_{2m+1}(-1/\kappa)\right).
\end{align}
Thus, if we have information about the arithmetic nature of $F_{2m+1}(\kappa)-\kappa^{2m}F_{2m+1}(-1/\kappa)$, that would help shed light on the arithmetic nature of the odd zeta values. In this direction, there is a result of Gun, Murty and Rath \cite{gmr} which says that if $m\in\mathbb{N}\cup\{0\}$ and $\delta_m=0, 1, 2$ or $3$ respectively as $\textup{gcd}(m, 6)=1, 2, 3$ or $6$, then for every algebraic $\kappa\in\mathbb{H}$, the number  $F_{2m+1}(\kappa)-\kappa^{2m}F_{2m+1}(-1/\kappa)$ is transcendental with at most $2m+2+\delta_m$ exceptions.

O'Sullivan \cite[Equation (8.24)]{sullivan} has obtained another representation for $\zeta(2m+1)$ in terms of the Ramanujan polynomials using his non-holomorphic analogue of \eqref{rameqn}, that is, Theorem \ref{non-holomorphic-thm}. It states that for $m\in\mathbb{N}$,
\begin{align}
	\zeta(2m+1)&=\frac{1}{(z^{2m}-1)}\bigg\{\frac{1}{z}(2\pi i)^{2m+1}R_{2m+1}(z)+2\left(V_{-2m}(z)-z^{2m}V_{-2m}(-1/z)\right)\nonumber\\
	&\quad+\frac{B_{2m+2}}{2(2m+2)!}z^{2m}(|z|^{-4m-2}-1)\bigg\}.
\end{align}
If we let $\alpha=-i\pi z$ and $\beta=i\pi/z$ in each of the finite sums occurring in \eqref{equiforma=-2mano} and \eqref{trans2m+1eqn}, we get cousins of Ramanujan polynomials \eqref{rp}. It may be worthwhile to investigate their properties as well.\\

\section{Conclusion}

One of the objectives of this paper was to demonstrate how rich the topic concerning Ramanujan's formula for odd zeta values is. There have been incredibly many generalizations and ramifications of it, and while newer perspectives continue to emerge, sometimes older results also get rediscovered.   

For example, consider the following equivalent representation of \eqref{rameqn} given in \cite[Equation (1)]{uhl}\footnote{This rephrasing of \eqref{rameqn} was done first by Ramanujan himself. See also \cite[p.~155]{berndtrocky}.}:
\begin{align}\label{uhl1}
	\alpha^{m-1}\sum_{n=1}^{\infty}\frac{\coth(\pi n\alpha^{-1})}{(2\pi n)^{2m-1}}	-(-\alpha)^{1-m}\sum_{n=1}^{\infty}\frac{\coth(\pi n\alpha)}{(2\pi n)^{2m-1}}=-\sum_{n=0}^{m}(-1)^n\frac{B_{2n}}{(2n)!}\frac{B_{2m-2n}}{(2m-2n)!}\alpha^{2n-m},
\end{align}
where we consider $m\geq2$. Define
\begin{align}
	C_m(\alpha):=\sum_{n=1}^{\infty}\frac{\coth(\pi n\alpha)}{(2\pi n)^{2m-1}},\hspace{6mm} G_m(\alpha):=-\sum_{n=0}^{m}(-1)^n\frac{B_{2n}}{(2n)!}\frac{B_{2m-2n}}{(2m-2n)!}\alpha^{2n-m}.
	\end{align}
Then \eqref{uhl1} implies
\begin{equation}\label{rameqn analogue}
	G_m(\alpha)=\alpha^{m-1}C_{m}(\alpha^{-1})-(-\alpha)^{1-m}C_m(\alpha).
\end{equation}
Uhl gives a new proof of \eqref{rameqn analogue} using Mittag-Leffler expansion\footnote{We note that Chavan \cite{schavan} has also obtained a new proof of \eqref{rameqn} using Mittag-Leffler expansion.} and then obtains \cite[Theorem 4.1]{uhl}, what he calls, a triangle identity for $F_m(\alpha)$:
\begin{equation*}
	G_m(\alpha)=(i\alpha^{-1}+1)^{m-1}G_m(\alpha+i)-(i\alpha-1)^{m-1}G_m(\alpha^{-1}-i).
\end{equation*}
It is not difficult to see that the above identity was already discovered in an equivalent form by Vlasenko and Zagier \cite[p.~42, Equation (30)]{vz}.

In the abstract of his paper, Uhl writes, `\emph{\dots properties and symmetries of the equation [Ramanujan's formula] are far from all uncovered}'.

 While the scope of the current survey was quite limited, the above quote certainly hints at a new survey article to emerge in the years to come! 

\begin{center}
	\textbf{Acknowledgements}
\end{center}
We sincerely thank Bruce C. Berndt, Rahul Kumar and N. Guru Sharan for carefully going over the manuscript and for giving nice suggestions.


\begin{thebibliography}{00}
	\bibitem{as}
	M.~Abramowitz and I.A.~Stegun, eds., \emph{Handbook of
		Mathematical Functions}, Dover, New York, 1965.
	
		\bibitem{RLNII}
	George E.~Andrews and Bruce C.~Berndt, \emph{Ramanujan's Lost Notebook Part \textup{II}}, Springer, New York, 2009, 418 pp.
	
	\bibitem{apery1}
	R.~Ap\'{ery}, \emph{Irrationalit\'{e} de $\zeta(2)$ et $\zeta(3)$}, Ast\'{e}risque~\textbf{61} (1979), 11--13.
	
	\bibitem{apery2}
	R.~Ap\'{e}ry, \emph{Interpolation de fractions continues et irrationalit\'{e} de certaines constantes},
	Bull. Section des Sci., Tome III, Biblioth\'{e}que Nationale, Paris, 1981, 37--63.
	
	\bibitem{apostol-1998a}
	Tom M.~Apostol, \emph{Introduction to Analytic Number Theory}, Springer-Verlag, New York (1998).
	
	\bibitem{apostol2}
	Tom M.~Apostol, \emph{Modular Functions and Dirichlet Series in Number Theory}, $2^{\textup{nd}}$ edition, Springer, 1990.
	
	\bibitem{ballrivoal}
	K.~Ball and T.~Rivoal, \emph{Irrationalit\'{e} d'une infinit\'{e} de valeurs de la fonction z\^{e}ta aux entiers impairs} (French), Invent. Math.~\textbf{146} no. 1 (2001), 193--207.
	
		\bibitem{bdg_log}
	S.~Banerjee, A.~Dixit and S.~Gupta, \emph{Lambert series of logarithm, the derivative of Deninger's function $R(z)$ and a mean value theorem for $\zeta\left(\frac{1}{2}-it\right)\zeta'\left(\frac{1}{2}+it\right)$}, Canad.~J.~Math.~(2023) (36 pages) (DOI: \url{https://doi.org/10.4153/S0008414X23000597})
	
	\bibitem{BDG}
	S.~Banerjee, A.~Dixit and S.~Gupta, \emph{Explicit transformations for generalized Lambert series associated with the divisor function $\sigma_a^{(N)}(n)$ and their applications}, Res.~Math.~Sci.~\textbf{10} no. 4, (2023),  Paper No. 38 (50 pages).
	
		\bibitem{bgk}
	S.~Banerjee, R.~Gupta and R.~Kumar, \emph{A note on odd zeta values over any number field and Extended Eisenstein series}, J.~Math.~Anal.~Appl.~\textbf{531} (2024), 127883 (17 pp.)
	
	\bibitem{banerjee-sahani}
	S.~Banerjee and V.~Sahani, \emph{Transformation formulas for the higher power of odd zeta values and generalized Eisenstein series}, submitted for publication \url{https://arxiv.org/pdf/2206.13331.pdf}.

	\bibitem{berndthurwitzzeta}
B.~C.~Berndt, \emph{On the Hurwitz zeta function}, Rocky~Mountain~J.~Math.~\textbf{2} (1972), 151--157.

\bibitem{berndtrocky}
	B.~C.~Berndt, \emph{Modular transformations and generalizations of several formulae of Ramanujan}, Rocky Mountain J. Math.~\textbf{7} (1977), 147--189.

\bibitem{berndt-jupiter}	
	B.~C.~Berndt, \emph{Ramanujan's formula for $\zeta (2n+1)$}, in
	{\it Professor Srinivasa Ramanujan
		Commemoration Volume}, Jupiter Press, Madras, 1974, pp.~1--7.
		
	\bibitem{berndtcrelle}
B.~C.~Berndt, \emph{Analytic Eisenstein series, theta-functions, and series relations in the spirit of Ramanujan}, 
J.~Reine~Angew.~Math.~\textbf{303}(304) (1978), 332--365.
	
	\bibitem{bcbramsecnote}
	B.~C.~Berndt, \emph{Ramanujan's Notebooks, Part II}, Springer-Verlag, New York, 1989.
	
	\bibitem{bcbramthinote}
	B.~C.~Berndt, \emph{Ramanujan's Notebooks, Part III}, Springer-Verlag, New York, 1991.
	
	\bibitem{bcbramfifthnote}
	B.~C.~Berndt, \emph{Ramanujan's Notebooks, Part V}, Springer-Verlag, New York, 1998.
	
	
%	\bibitem{bct}
%	B.~C.~Berndt, H.~H.~Chan and Y.~Tanigawa, \emph{Two Dirichlet series evaluations found on page 196 of Ramanujan's Lost Notebook}, Math. Proc. Camb. Phil. Soc.~\textbf{153} No. 2 (2012), 341--360.
	
		\bibitem{bdgz}
	B.~C.~Berndt, A.~Dixit, R.~Gupta and A.~Zaharescu, \emph{Ramanujan and Koshliakov meet Abel and Plana}, In: Baskar Balasubramanyam, Kaneenika Sinha and Mathukumalli Vidyasagar (eds.) \emph{Some contributions to number theory and
	beyond: Proceedings of the centenary symposium for M V Subbarao},  Fields Institute Communications (to appear).
%	
%	\bibitem{bdrz1}
%	B.~C.~Berndt, A.~Dixit, A.~Roy and A.~Zaharescu, \emph{New pathways and connections in number theory and analysis motivated by two incorrect claims of Ramanujan}, Adv. Math.~\textbf{304} (2017), 809--929.
	
	\bibitem{berndtstraubmathz}
	B.~C.~Berndt and A.~Straub, \emph{On a secant Dirichlet series and Eichler integrals of Eisenstein series}, Math. Z.~\textbf{284} No. 3-4 (2016), 827--852.
	
	\bibitem{berndtstraubzeta}
	B.~C.~Berndt and A.~Straub, \emph{Ramanujan's formula for $\zeta(2n+1)$}, Exploring the Riemann zeta function, Eds. H. Montgomery, A. Nikeghbali, and M. Rassias, pp.~13--34, Springer, 2017.
	
%	\bibitem{boxalljones}
%	G.~J.~Boxall and G.~O.~Jones, \emph{Algebraic values of certain analytic functions}, Int. Math. Res. Not.~ IMRN (2015), No. 4, 1141--1158.
%	
%	\bibitem{hmfmmf}
%	K.~Bringmann, A.~Folsom, K.~Ono, and L.~Rolen, \emph{Harmonic Maass Forms and
%		Mock Modular Forms: Theory and Applications}, Amer. Math. Soc., Providence, 2017.
%	
%	\bibitem{cassels}
%	J.~W.~S.~Cassels, \emph{Footnote to a note of Davenport and Heilbronn}, J. London Math. Soc.~\textbf{36} (1961), 177--184.
	
	\bibitem{schavan}
	S.~Chavan, \emph{An elementary proof of Ramanujan’s identity for odd zeta values}, J.~Class.~Anal.~\textbf{19} (2) (2022), 139--147.
	
		\bibitem{pchavan}
	P.~Chavan, \emph{Hurwitz zeta functions and Ramanujan's identity for odd zeta values}. J.~Math.~Anal.~Appl.~\textbf{527} no. 2 (2023), 20 pp.
	
	\bibitem{ccvw}
	P.~Chavan, S.~Chavan, C.~Vignat and T.~Wakhare, \emph{Dirichlet series under standard convolutions: variations on Ramanujan's identity for odd zeta values}, Ramanujan J.~\textbf{59} no. 4 (2022), 1245--1285.
			
	\bibitem{cjm}
	S.~Chourasiya, Md. K.~Jamal and B.~Maji, \emph{A new Ramanujan-type identity for $L(2k+1,\chi_1)$}, Ramanujan~J.~\textbf{60} no. 3 (2023), 729--750. 
		
		\bibitem{conrey}
		J.~B.~Conrey, D.~W.~Farmer and  \"{O}.~Imamoglu, \emph{The nontrivial zeros of period polynomials of modular forms lie on the unit circle}, Int. Math. Res. Not. IMRN\textbf{20} (2013), 4758--4771.
%	\bibitem{dainaylor}
%	H.~H.~Dai and D.~Naylor, \emph{On an asymptotic expansion of Fourier integrals}, Proc. Roy. Soc. London Ser. A~\textbf{436} No. 1896 (1992), 109--120.

\bibitem{dagli-can}
M.~C.~Da\u{g}li and M.~Can, \emph{On generalized Eisenstein series and Ramanujan's formula for periodic zeta-functions}, Monatsh.~Math.~\textbf{184} no. 1 (2017), 77--103.

	\bibitem{danbiyarkhodzhaev}
	A.~T.~Daniyakhodzhaev and M.~A.~Korolev, \emph{On a Ramanujan Identity and Its Generalizations }, Mat.~Zmetki~\textbf{110} no. 4 (2021), 524--536; Math. Notes~\textbf{110} no. 4 (2021), 511--521.
%	
%	\bibitem{dav}
%	H.~Davenport, \emph{Multiplicative Number Theory}, 3rd ed., Springer--Verlag, New York, 2000.
%	
%	\bibitem{davhel}
%	H.~Davenport and H.~Heilbronn, \emph{On the zeros of certain Dirichlet series, I}, J. London Math. Soc.~\textbf{11} (1936), 181--185.
	
	\bibitem{dilcher}
	K.~Dilcher, \emph{On generalized gamma functions related to the Laurent coefficients of the Riemann zeta function}, Aequationes Math.~\textbf{48} (1994), 55--85.
%	
%	\bibitem{dixitlap}
%	A.~Dixit, \emph{The Laplace transform of the psi function}, Proc. Amer. Math. Soc.~\textbf{138} No. 2 (2010), 593--603.
			
		
		\bibitem{dgsquare}
		A. Dixit and R. Gupta, \emph{On squares of odd zeta values and analogues of Eisenstein series}, Adv.~Appl.~Math.~\text{110} (2019), 86--119.
			
			\bibitem{dg}
			A. Dixit and R. Gupta, \emph{Koshliakov zeta functions I. Modular relations}, Adv. Math. \textbf{393} (2021), Paper No. 108093 (41 pages).
			
			\bibitem{dgk}
			A.~Dixit, R.~Gupta and R.~Kumar, \emph{Extended higher Herglotz functions I. Functional equations}, Adv. ~Appl.~Math.~\textbf{153} (2024), 102622 (41 pages).
			
			\bibitem{dgkm}
			A.~Dixit, R.~Gupta, R.~Kumar and B.~Maji, \emph{Generalized Lambert series, Raabe's cosine transform and a generalization of Ramanujan's formula for $\zeta(2m+1)$}, Nagoya Math.~J.~\textbf{239} (2020), 232--293.
			
	\bibitem{dkk}
	A.~Dixit, A.~Kesarwani and R.~Kumar, \emph{Explicit transformations of certain Lambert series}, Res.~ Math.~Sci.~\textbf{9}, 34 (2022) (54 pages). 
	
		\bibitem{dk03}
	A.~Dixit and R.~Kumar, \emph{Applications of the Lipschitz summation formula and a generalization of Raabe's cosine transform}, Constr. Approx. (2023) (40 pages). (DOI: \url{https://doi.org/10.1007/s00365-023-09668-8})
	
\bibitem{dixitmaji1}
A.~Dixit and B.~Maji, \emph{Generalized Lambert series and arithmetic nature of odd zeta values}, Proc.~Royal~Soc.~Edinburgh, Sect. A: Mathematics, \textbf{150} Issue 2 (2020), 741--769.
	
	\bibitem{dgs1}
	A.~Dixit, N.~G.~Sharan and S.~Sathyanarayana, \emph{Modular relations involving generalized digamma functions}, submitted for publication. \url{https://arxiv.org/pdf/2306.10991.pdf}
	
	\bibitem{dorigoni}
	D.~Dorigoni and A.~Kleinschmidt, \emph{Resurgent expansion of Lambert series and iterated Eisenstein integrals}, Commun. Number Theory Phys.~\textbf{15} No. 1 (2021), 1--57.
	
%	\bibitem{dixfer1}
%	A.~L.~Dixon and W.~L.~Ferrar, \emph{Lattice-point summation formulae}, Quart.~J.~Math.~\textbf{2} (1931), 31--54.
%	
%	\bibitem{dixfer3}
%	A.~L.~Dixon and W.~L.~Ferrar, \emph{Infinite integrals of Bessel functions}, Quart.~J.~Math.~\textbf{1} (1935), 161--174.
	
%\bibitem{eie-chen}
%M.~Eie and K.-W.~Chen, \emph{A theorem on zeta functions associated with polynomials}, Trans.~Amer.~Math.~Soc.~\textbf{351} no. 8 (1999), 3217--3228.

%	\bibitem{htf2}
%	A.~Erd\'{e}lyi, W.~Magnus, F.~Oberhettinger and F.~Tricomi, \emph{Higher Transcendental Functions} (Bateman Manuscript Project), Vol. II, McGraw-Hill, New York, 1953.
%	
%	\bibitem{erdos}
%	P.~Erd\"{o}s, \emph{On arithmetical properties of Lambert series}, J. Indian Math. Soc. (N. S.) \textbf{12} (1948), 63--66.
	\bibitem{transseries}
	G.~A.~Edgar, \emph{Transseries for beginners}, Real Anal.~Exchange~\textbf{35} no. 2 (2010), 253--309.
	
	\bibitem{franke-rj}
		J.~Franke, \emph{Infinite series representations for Dirichlet L-functions at rational arguments}, Ramanujan J.~\textbf{46} no. 1 (2018), 91--102.
	
	\bibitem{franke-rint}
	J.~Franke, \emph{Ramanujan identities of higher degree}, Res.~Number Theory~\textbf{4} no. 4 (2018), 42, 19 pp.
	
	\bibitem{guptakumar}
	R.~Gupta and R.~Kumar, \emph{Extended higher Herglotz function II}, J.~Math.~Anal.~Appl.~\textbf{518} (2023), 12670 (16 pp.)
		
		\bibitem{gjkm}
		A.~Gupta, Md.~K.~Jamal, N.~Karak and B.~Maji, \emph{A Dirichlet character analogue of Ramanujan's formula for odd zeta values}, submitted for publication \url{https://arxiv.org/pdf/2308.08988.pdf}.
		
			\bibitem{gmr}
		S.~Gun, M.~R.~Murty and P.~Rath, \emph{Transcendental values of certain Eichler integrals}, Bull. London Math. Soc.~\textbf{43} No. 5 (2011), 939--952.
		
		\bibitem{guptamaji}
	A.~Gupta and B.~Maji, \emph{On Ramanujan’s formula for $\zeta(1/2)$ and $\zeta(2m+1)$}, J.~Math.~Anal.~Appl.~\textbf{507} (2022), 125738.
	
	\bibitem{grn}
	I.~S.~Gradshteyn and I.~M.~Ryzhik, eds., \emph{Table of Integrals,
		Series, and Products}, 8th ed., Edited by D.~Zwillinger, V.~H.~Moll, Academic Press, New York, 2015.
	
	\bibitem{gross1}
	E.~Grosswald, \emph{Die Werte der Riemannschen Zetafunktion an ungeraden Argumentstellen}, Nachr.~Akad.~Wiss.~G\"{o}ttinger Math.-Phys. Kl. II (1970), 9--13.
	
	\bibitem{gross2}
	E.~Grosswald, \emph{Comments on some formulae of Ramanujan}, Acta Arith.~\textbf{21} (1972), 25--34.
	
%	\bibitem{gross3}
%	E.~Grosswald, \emph{Relations between the values at integral arguments of Dirichlet series that satisfy functional equations}, Proc. Sympos. Pure Math., Vol. 24, Amer. Math. Soc., Providence, 1973, 111--122.
	
%	\bibitem{apg1}
%	A.P.~Guinand, \emph{On Poisson's summation formula}, Ann.~Math.
%	(2) \textbf{42} (1941), 591--603.
%	
%	\bibitem{gonek}
%	S.~M.~Gonek, \emph{Analytic Properties of zeta and $L$-functions}, Thesis, University of Michigan, 1979.


%	\bibitem{hanclkristensen}
%	J.~Han\v{c}l and S.~Kristensen, \emph{Metrical irrationality results related to values of the Riemann $\zeta$-function}, arXiv:1802.03946v1, February 12, 2018.
	
%	\bibitem{je}
%	E.~Jahnke and F.~Emde, \emph{Tables of functions with formulae and curves}, 4th ed., Dover Publications, New York, 1945.
	
	\bibitem{jmos}
	S.~Jin, W.~Ma, K.~Ono and K.~Soundararajan, \emph{Riemann hypothesis for period polynomials of modular forms}, Proc.~Nat.~Acad.~Sci~USA~\textbf{113} no. 10 (2016), 2603--2608.
	
	\bibitem{ktyhr}
	S.~Kanemitsu, Y.~Tanigawa and M.~Yoshimoto, \emph{On the values of the Riemann zeta-function at rational arguments}, Hardy-Ramanujan J.~\textbf{24} (2001), 11--19.
	
	\bibitem{ktyham}
	S.~Kanemitsu, Y.~Tanigawa, and M.~Yoshimoto, \emph{On rapidly convergent series for the Riemann zeta-values via the modular relation}, Abh. Math. Sem. Univ. Hamburg~\textbf{72} (2002), 187--206.
	
	\bibitem{ktyacta}
	S.~Kanemitsu, Y.~Tanigawa, and M.~Yoshimoto, \emph{On multiple Hurwitz zeta-function values at rational arguments}, Acta Arith.~\textbf{107}, No. 1 (2003), 45--67.
	
	\bibitem{katsurada}
	M.~Katsurada and T.~Noda, \emph{Asymptotic expansions for a class of generalized holomorphic Eisenstein series, Ramanujan's formula for $\zeta(2k+1)$, Weierstrass' elliptic and allied functions}, submitted for publication \url{https://arxiv.org/pdf/2201.10124.pdf}
	
	\bibitem{kirprod}
	P.~Kirschenhofer and H.~Prodinger, \emph{On some applications of formulae of Ramanujan in the analysis of algorithms}, Mathematika~\textbf{38} No. 1 (1991), 14--33.
	
	\bibitem{kongsiriwong}
	S.~Kongsiriwong, \emph{A generalization of Siegel's method}, Ramanujan J.~\textbf{20} (2009), 1--24.
	
		\bibitem{koshliakov3}
	N.~S.~Koshliakov (under the name N.S.~Sergeev), \emph{Issledovanie odnogo klassa transtsendentnykh
		funktsii, opredelyaemykh obobshcennym yravneniem Rimana} (A study of a
	class of transcendental functions defined by the generalized
	Riemann equation) (in Russian), Trudy Mat. Inst. Steklov, Moscow,
	1949. (Available online at \url{https://dds.crl.edu/crldelivery/14052})
	
%	\bibitem{koyakuro}
%	S.~Koyama and N.~Kurokawa, \emph{Kummer's formula for multiple gamma functions}, J.~Ramanujan Math. Soc.~\textbf{18} No. 1 (2003), 87--107.
%	
%	\bibitem{kummer}
%	E.~E.~Kummer, \emph{Beitrag zur Theorie der Function $\G(x)=\int_{0}^{\infty}e^{-v}v^{x-1}dv$}, J. Reine Angew. Math.~\textbf{35} (1847), 1--4.
	
%	\bibitem{lagariaseuler}
%	J.~Lagarias, \emph{Euler's constant: Euler's work and modern developments}, Bull. Amer. Math. Soc.~\textbf{50} No. 4 (2013), 527--628.
%	
%	\bibitem{lagrange}
%	J.~Lagrange, \emph{Une formule sommatoire et ses applications}, Bull.~Sci.~Math.~(2)~\textbf{84} (1960), 105--110.
	
%	\bibitem{lerchform}
%	M.~Lerch, \emph{Dal\v{s}i studie v oboru Malmst\'{e}novsk\'{y}ch \v{r}ad}, Rozpravy \v{C}esk\'{e} Akad.~\textbf{3} (28) (1894), 1--61.
%	
\bibitem{lalin-rogers}
M. N. Lal\'{i}n and M. D. Rogers, \emph{Variations of the Ramanujan polynomials and remarks on $\zeta(2j+1)/\pi^{2j+1}$}, Funct.~Approx.~Comment.~Math.~\textbf{48}, no. 1 (2013), 91--111.

\bibitem{lalin-smyth}
M. N. Lal\'{i}n and C.~J.~Smyth, \emph{Unimodularity of zeros of self-inversive polynomials}, Acta~Math.~Hungar.~\textbf{138} no. 1-2 (2013), 85--101; Addendum, Acta Math. Hungar.~\textup{147} no. 1 (2015), 255--257.

	\bibitem{lerch}
	M.~Lerch, \emph{Sur la fonction $\zeta(s)$ pour valeurs impaires de l'argument}, J. Sci. Math. Astron. pub. pelo Dr. F. Gomes Teixeira, Coimbra~\textbf{14} (1901), 65--69.
	
%	\bibitem{lucatachiya}
%	F.~Luca and Y.~Tachiya, \emph{Linear independence of certain Lambert series}, Proc. Amer. Math. Soc.~\textbf{142} No. 10 (2014), 	3411--3419.
%	

%\bibitem{lim-acta}
%S.-G.~Lim, \emph{Infinite series Identities from modular transformation formulas that stem from generalized Eisenstein series}, Acta~Arith.~\textbf{141} (3) (2010), 241--273.

\bibitem{lim1}
S.-G.~Lim, \emph{On some series identities}, Honam Math.~J.~\textbf{38} no. 3 (2016), 479--494.

	\bibitem{malurkar}
	S.~L.~Malurkar, \emph{On the application of Herr Mellin’s integrals to some series}, J.~Indian~Math.~Soc.~\textbf{16}, (1925/26), 130--138.
	
	\bibitem{maier}
	W.~Maier, \emph{Abbau singul\"{a}rer linien}, Math.~Ann.~\textbf{169} (1967), 102--117.
	
%	\bibitem{magobersoni}
%	W.~Magnus, F.~Oberhettinger and R.~P.~Soni, \emph{Formulas and Theorems for the Special Functions of Mathematical Physics}, 3rd ed., Vol. 52, Springer-Verlag, New York, 1966.
	
%	\bibitem{masser}
%	D.~Masser, \emph{Rational values of the Riemann zeta function}, J. Number Theory~\textbf{131} (2011), 2037--2046.
	
	\bibitem{msw}
	M.~R.~Murty, C.~Smyth and R.~J.~Wang, \emph{Zeros of Ramanujan polynomials}, J.~Ramanujan Math.~Soc.~\textbf{26} No. 1 (2011), 107--125.
	
%	\bibitem{ober}
%	F.~Oberhettinger, \emph{Tables of Mellin Transforms}, Springer-Verlag, New York, 1974.
	
%	\bibitem{olver1974}
%	F.~W.~J.~Olver, \emph{Error bounds for stationary phase approximations}, SIAM J. Math. Anal.~\textbf{5} No. 1 (1974), 19--29.
	
	\bibitem{sullivan}
	C.~O'Sullivan, \emph{Formulas for non-holomorphic Eisenstein series and for the Riemann zeta function at odd integers}, Res.~Number Theory~\textbf{4} no. 3 (2018), Art. 36, 38 pp.
	\bibitem{olver-2010a}
	F.~W.~J.~Olver, D.~W.~Lozier, R.~F.~Boisvert, and C.~W.~Clark, eds., \emph{NIST Handbook of Mathematical Functions}, Cambridge University Press, Cambridge, 2010.
	
%	\bibitem{kp}
%	R.~B.~Paris and D.~Kaminski, \emph{Asymptotics and Mellin-Barnes Integrals},  Encyclopedia of Mathematics and its Applications, 85. Cambridge University Press, Cambridge, 2001.
%	
%	\bibitem{prudv3}
%	A.~P.~Prudnikov, Yu.~A.~Brychkov and O.~I.~Marichev, \emph{Integrals and
%		Series, Vol.~3, More Special Functions}, Gordon and Breach, New York, 1986.
%	
%	\bibitem{rademacheraddress}
%	H.~Rademacher, \emph{Trends in research: The Analytic Number Theory}, Address delivered by invitation of the American Mathematical Society Program Committee, September 5, 1941; Bull. Amer. Math. Soc.~\textbf{48} (1942), 379--401.
	
	\bibitem{rajkumar}
	K.~Rajkumar, \emph{A simplification of Ap\'{e}ry’s proof of the irrationality of $\zeta(3)$}, (2012), 7 pp. \url{https://arxiv.org/pdf/1212.5881.pdf}
	
	\bibitem{ramnote}
	S.~Ramanujan, Notebooks (2 volumes), Tata Institute of Fundamental Research, Bombay, 1957; second ed., 2012.
	
	\bibitem{lnb}
	S.~Ramanujan, \emph{The Lost Notebook and Other Unpublished
		Papers}, Narosa, New Delhi, 1988.
	
	\bibitem{rivoal}
	T.~Rivoal, \emph{La fonction z\^{e}ta de Riemann prend une infinit\'{e} de valeurs irrationnelles aux entiers impairs}, C. R. Acad. Sci. Paris S\'{e}r. I Math.~\textbf{331} no. 4 (2000), 267--270.
	
%	\bibitem{schlomilch}
%	O.~Schl\"{o}milch, \emph{Ueber einige unendliche Reihen}, Berichte \"{u}ber die Verh. d. K\"{o}nige S\"{a}chsischen Gesell. Wiss. zu Leipzig~\textbf{29} (1877), 101--105.
%	
%	\bibitem{sitaramachandrarao}
%	R.~Sitaramachandrarao, \emph{Ramanujan’s formula for $\zeta(2n+1)$}, Madurai Kamaraj University
%	Technical Report 4, pp. 70--117.
%	
%	\bibitem{spira}
%	R.~Spira, \emph{Zeros of Hurwitz zeta functions}, Math.~Comp.~\textbf{30} No. 136 (1976), 863--866.
	
%	\bibitem{srivastava}
%	H.~M.~Srivastava, \emph{Further series representations for $\zeta(2n+1)$}, Appl. Math. Comput.~\textbf{97} (1998), 1--15.
%	

\bibitem{straub}
A.~Straub, \emph{Special values of trigonometric Dirichlet series and Eichler integrals}, Ramanujan J.~\textbf{41}, no. 1-3 (2016), 269--285.
 
\bibitem{temme}
	N.~M.~Temme, \emph{Special functions: An introduction to the classical functions of mathematical physics}, Wiley-Interscience Publication, New York, 1996.

%	\bibitem{terras0}
%	A.~Terras, \emph{Some formulas for the Riemann zeta function at odd integer argument resulting from Fourier expansions of the Epstein zeta function}, Acta Arith.~\textbf{29} No. 2 (1976), 181--189.
%	
%	\bibitem{terras}
%	A.~Terras, \emph{The Fourier expansion of Epstein's zeta function over an algebraic number field and its consequences for algebraic number theory}, Acta Arith.~\textbf{32} (1977), 37--53.
%	
%	\bibitem{titchfou}
%	E.~C.~Titchmarsh, \emph{Theory of Fourier Integrals}, 2nd ed., Clarendon Press, Oxford, 1948.
	
	\bibitem{uhl}
	M.~Uhl, \emph{Ramanujan’s formula for odd zeta values: a proof by Mittag-Leffler expansion and applications}, Eur.~J.~Math.~\textbf{9} no. 3 (2023), Paper No. 79, 12 pp.
	
	\bibitem{vz}
	M.~Vlasenko and D.~Zagier, \emph{Higher Kronecker ``limit" formulas for real quadratic fields}, J. reine angew. Math. 679, pp. 23--64 (2013).
	
%	\bibitem{voronin}
%	S.~M.~Voronin, \emph{On the zeros of zeta-functions of quadratic forms}, Trudy Mat. Inst. Steklov~\textbf{142} (1976), 135--147; English translation in Proc. Steklov Inst. Math.~\textbf{3} (1979), 143-155.
%	
%	\bibitem{waldschmidt}
%	M.~Waldschmidt, \emph{Transcendence of periods: the state of the art}, Pure Appl. Math. Q.~\textbf{2} (2) (2006), 435--463.
	
	\bibitem{wig0}
	S.~Wigert, \emph{Sur la s\'{e}rie de Lambert et son application $\grave{a}$ la th\'{e}orie des nombres}, Acta Math.~\textbf{41} (1916), 197--218.
	
%	\bibitem{wig}
%	S.~Wigert, \emph{Sur une extension de la s\'{e}rie de Lambert}, Arkiv Mat.~Astron.~Fys.~\textbf{19} (1925), 13 pp.
	
%	\bibitem{wiltonmess}
%	J.~R.~Wilton, \emph{A proof of Burnside's formula for $\log\G(x+1)$ and certain allied properties of Riemann's $\zeta$-function}, Mess. Math.~\textbf{52} (1922/1923), 90--93.
	
	\bibitem{zudilin}
	W.~W.~Zudilin, \emph{One of the numbers $\zeta(5), \zeta(7), \zeta(9)$ and $\zeta(11)$ is irrational} (Russian), Uspekhi Mat. Nauk~\textbf{56} No. 4 (2001), 149--150; translation in Russian Math. Surveys~\textbf{56} No. 4 (2001), 774--776.
	
	\bibitem{zudilin1}
	W.~Zudilin, \emph{Ramanujan and odd zeta values}, to appear in the \emph{Encyclopedia of Srinivasa Ramanujan and His Mathematics} \url{https://wain.mi-ras.ru/PS/ramaodd.pdf}
\end{thebibliography}
\end{document}